\documentclass[10pt]{article}
\usepackage{oldgerm}
\usepackage{latexsym}
\usepackage{amsmath}
\usepackage{amssymb}
\usepackage{amsgen}
\usepackage{amscd} 

\input{diagrams} 
\diagramstyle[h=7mm,PostScript=dvips]

\newif\ifextendedversion
\extendedversionfalse

\newtheorem{prop}{Proposition}[section]
\newtheorem{definition}[prop]{Definition}
\newtheorem{cor}[prop]{Corollary}
\newtheorem{lemma}[prop]{Lemma}
\newtheorem{thm}[prop]{Theorem}
\newcommand{\opr}[1]{\operatorname{#1}}
\newcommand{\transfer}{\opr{\stackrel{\star}{\phantom{.}}}}
\newcommand{\MDLM}{V(\mathbb{R})}

\newcommand{\ESS}{S}
\newcommand{\SCALE}{\sigma}
\newcommand{\Rescale}{\gamma}

\newcommand{\fraco}[2]{#1 / #2}

\newtheorem{ExampleAux}{Example}[section]

\newenvironment{example}{\begin{ExampleAux} \rm} {\end{ExampleAux}}

\newcommand{\phull}[1]{\opr{\diamond}(#1)} 
\newcommand{\lhull}[1]{\opr{\diamond}#1 }

\newarrow{Line} -----

\title{Nonstandard Analysis of Graphs} \author{F. Javier Thayer}

\begin{document} \maketitle \begin{abstract}
This paper shows certain classes of metric length spaces characterized by
volume growth properties of balls can viewed as graphs with
infinitesimal edges. Our approach is based on nonstandard analysis.

\end{abstract}

\section{Introduction}

Intuitively, it is tempting to view a cube in $\mathbb{R}^n$ as a
discrete lattice with infinitesimal edges. A similar infinitesimal
picture suggests itself with other fractal spaces such as the
Sierpinski gasket or the Sierpinski sponge~(see
\cite{morgan},~\cite{falconer}.)  This paper shows certain classes of
metric spaces characterized by a polynomial growth condition on the
volume of balls~(\cite{davidsemmes},\cite{chavel}) can be viewed as a
part of a graph with infinitesimal edges.  Our approach uses
nonstandard analysis in two ways: first, nonstandard analysis makes it
possible to do infinitesimal rescaling of graph distances in a
completely straightforward way, and second, we use the well-known
connection between nonstandard counting measures and countably
additive measures established in~\cite{Loeb75} and~\cite{henson} to
associate Hausdorff measure on metric spaces to counting measures on
graphs.  The main result is Theorem~\ref{main-characterization-result}
which exhibits an Ahlfors regular length space up to Lipschitz
equivalence as a part of a hyperfinite graph in which the edges have
infinitesimal length.

The structure of the paper is as follows: the next two sections
consist mainly of reference material, particularly on nonstandard
analysis and a summary of results in nonstandard measure theory.  In
particular, in Section~2 we give a definition of polynomial growth for
functions defined on intervals with hyperreal endpoints.  Note that
for functions which are defined on an interval whose endpoints are
real numbers, polynomial growth reduces to boundedness on that
interval.

In Section~\ref{polynomial-growth-section}, we introduce the main
concept of the paper, namely that of polynomial growth for hyperfinite
graphs.  The basic idea is as follows: Consider the shortest-path
metric between nodes.  If the cardinality of balls as a function of
the radius behaves polynomially in some interval, we say the graph has
polynomial growth on that interval. Note the notion of polynomial
growth for finite graphs is vacuous.  Along with the shortest-distance
metric, we consider rescalings of that metric and using nonstandard
analysis we have a lot of leeway on how we choose the rescaling
parameter.  In particular, if we choose an infinitesimal rescaling, we
get a graph in which nodes may be at an infinitesimal distance from
each other. Now collapse the graph, identifying nodes that are
infinitely close to each other.  The resulting object is a metric
space, which in the polynomial growth case has a volume growth
property called Ahlfors regularity as shown in
Theorem~\ref{measure-characterization-polynomial-growth}.
The remainder of the paper shows that the converse (suitably stated)
also holds for the class of length spaces introduced by Gromov
in~\cite{gromov}.  

We note that the relation between Gromov's geometric ideas and
nonstandard analysis is not new.  Gromov himself uses ultraproducts in
a construction he calls the asymptotic cone over a metric space. An
explicit connection between nonstandard analysis and asymptotic
geometry has been established in the beautiful
paper~\cite{dries-wilkie}.  The authors of that paper used this
connection to provide a different proof of Gromov's theorem that a
finitely generated group of polynomial growth has a nilpotent subgroup
of finite index.  However, the relation between polynomial growth of
graphs and volume growth seems to be new.

\section{Preliminaries}

We will consider {\em extended metric spaces} in which the metric is
allowed to assume the value $+\infty$.  Note that this generalization
does not introduce any new topological behavior. Given an extended
metric space $(X,d)$, define an equivalence relation for pairs $x$,
$y$ of $X$ by $d(x,y) < \infty$.  We will call the equivalence classes
of this relation the {\em bounded components} of $(X,d)$.  The bounded
components of $(X,d)$ are open and closed sets.
\begin{definition} \label{length-space-definition}
An extended metric space $(X,d)$ is a {\em length space} iff for every
pair $x, y \in X$ with $d(x,y) < \infty$, there is an isometric map $f:
[0, d(x,y)] \rightarrow X$ with $f(0)=x, f(d(x,y)) = y$.
\end{definition}
This is a specialization of Gromov's definition (see~\cite{gromov},
D\'{e}finition 1.7), although by Th\'{e}or\`{e}me 1.10
of~\cite{gromov} the two definitions are equivalent for complete
locally compact metric spaces. For example, if $X$ is a geodesically
complete Riemannian manifold, by the Hopf-Rinow theorem~\cite{chavel},
Theorem 1.9, $(X,d)$ is a length space.
The following result is Th\'{e}or\`{e}me 1.10, (i) of~\cite{gromov}.
See also the first chapter of~\cite{gromov-tr}:
\begin{prop}\label{locally-compact-length-space-have-compact-balls}
If $(X,d)$ is a complete locally compact length space, all closed
balls $\overline{\opr{B}}(a,r)$ with $r < \infty$ are compact.
\end{prop}

\subsection{Ball Borel Structure} If $(X,d)$ is an extended metric
space, the smallest $\sigma$-algebra containing the open balls
$\opr{B}(a,r)$ with $r < \infty$ is the {\em ball Borel structure}.
If $(X,d)$ is separable, then since every open set in $X$ is a
countable union of open balls, the ball Borel structure is identical
to the usual Borel structure generated by the open sets.  In general,
the ball Borel structure has fewer sets than the Borel structure.  For
instance, if $(X,d)$ is an uncountable discrete metric space~(3.2.5
of~\cite{Dieudonne}), the only balls in $X$ are $\emptyset$, $X$ and
singletons.  Thus the ball Borel sets of $X$ are precisely those sets
which are countable or have countable complements whereas any subset
of $X$ is open and so any subset of $X$ is Borel.

If $(X,d)$ is a metric space, and $A \subseteq X$, then the ball Borel
structure on $(A,d)$ may not be the same as the ball Borel structure
of $(X,d)$ relativized to $A$.  For example, let $(A,d)$ be an
uncountable discrete metric space, $X = A \cup \{z\}$ where $z \not\in
A$ and $B \subseteq A$ such that $B$ and $A \setminus B$ are
uncountable. Extend $d$ to a metric $d'$ on $X$ by letting $d'(z, a)
=1$ if $a \in B$ and $d'(z, a) =2$ if $a \in A \setminus B$.
Since $B = \opr{\overline{B}}(z,1) \cap A$, $B$ is a member of the
relativized ball Borel structure, but is not a ball Borel set in $(A, d)$.

Despite this negative result, the relativization of the ball Borel
structure of $(X,d)$ to any bounded component $X_0$ is the ball Borel
structure of $(X_0,d)$: This follows from the remark that any ball
$\opr{B}(x,r)$ with $r < \infty$, is either a subset of $X_0$ or is
disjoint from $X_0$.

There is no reason to expect pleasant behavior from the ball Borel
structure for nonseparable metric spaces.  However, the metric spaces
of interest to us have the property that all bounded components are
separable or equivalently, that all open balls of finite radius are
separable.  For instance, it follows immediately from
Proposition~\ref{locally-compact-length-space-have-compact-balls},
that all bounded components of a length space are $\sigma$-compact.
In this case the relation between the ball Borel structure and the
Borel structure is easily determined:
\begin{prop}
Suppose $(X,d)$ is an extended metric space whose bounded components
are all separable.  Then the ball Borel structure of $(X,d)$ is the
$\sigma$-algebra $\mathcal{B}$ of those Borel sets $A$ such that $A$
or $\complement A$ is a subset of a countable union of bounded
components of $X$.
\end{prop}
{\sc Proof.}  $\mathcal{B}$ is clearly a $\sigma$-algebra.  Suppose
$A$ is a Borel set.  If $A$ is a subset of a bounded component $X_0$
of $X$, then by the equality of Borel structure and ball Borel
structure on separable spaces, $A$ is ball Borel in $X_0$ and
therefore ball Borel in $X$.  It follows that if $A$ or $\complement
A$ is a subset of a countable union of bounded components of $X$, then
$A$ is ball Borel.  Conversely, $\mathcal{B}$ contains all open balls
of finite radius, so contains all ball Borel sets.~$\blacksquare$

References for measure theory are~(\cite{dudley},~\cite{halmos}).  We
note one discrepancy between our notation and that used in these
references.  If $\mu$ is a countably additive measure on $(X,
\mathcal{A})$ and $g:(X, \mathcal{A}) \rightarrow (Y, \mathcal{B})$
is measurable then we use the notation $g_\ast \mu$ to denote the
measure $\mu \, g^{-1}$ on $(Y, \mathcal{B})$.

\subsection{Nonstandard Analysis}

We need a very small amount of background material on nonstandard
analysis such as the first few pages of Keisler's
monograph~\cite{keisler}.  Besides the Keisler reference, the
book~\cite{albeverio-etal} is also highly recommended and will be
referred to at various places in the paper.
Another general reference for this section with more foundational
material is
~\cite{hurd-loeb}.

The {\em superstructure} over $\mathbf{G}$ is the set $V(\mathbf{G})$
defined by:
$V_0(\mathbf{G}) = \mathbf{G}$, $V_{n+1}(\mathbf{G})
=V_{n}(\mathbf{G}) \cup \mathcal{P}(V_{n}(\mathbf{G}))$ and
$V(\mathbf{G}) = \bigcup_n V_{n}(\mathbf{G})$.
The main constituent of our working view of nonstandard analysis is a
map $\transfer: \MDLM \longrightarrow V(\transfer \mathbb{R})$ which
satisifies the transfer principle, c.f.~\cite{keisler} for details. In
this paper we will only use the countable saturation property.

If $\mathfrak{X}$ is an internal set, $\opr{card} \mathfrak{X}$
denotes the internal cardinality of $\mathfrak{X}$ if $\mathfrak{X}$
is hyperfinite, $+\infty$ otherwise.

If $r \in \transfer{\mathbb{R}}$, $r \ll \infty$ means that $r$ is
dominated by a standard real, $r \gg -\infty$ means that $r$ dominates
a standard real, $r$ is {\em limited} iff $-\infty \ll r \ll \infty$,
$r$ is {\em infinitesimal} iff for every positive standard real
$\theta$, $|r| \leq \theta$. Hyperreals $r, r'$ are {\em infinitely
close}, written $r \cong r'$ iff $r - r'$ is infinitesimal. $r \ll r'$
means $r' - r$ is positive and not infinitesimal.
The unique $r_0 \in \mathbb{R}$ infinitely close to $r$, if it exists is
the {\em standard part} of $r$ denoted $\opr{st}(r)$.
Define an (external) relation on $\transfer{[0, \infty[}$ by $r
\preceq_O r'$ iff there is a limited hyperreal $A$ such that $r \leq A
\ r'$.  Note that if $r' > 0$, then this is equivalent to
$\fraco{r}{r'} \ll \infty$.  Similarly, define $r \preceq_o r'$ iff
there is an infinitesimal hyperreal $A$ such that $r \leq A \ r'$.
Define $r \sim_O r'$ iff $r \preceq_O r'$ and $r' \preceq_O r$.  The
relation ``$\sim_O$'' is clearly an equivalence relation, albeit an
{\em external} one.  If $r \sim_O r'$ we say $r, r'$ are of the same
order of magnitude.
Note that for $r' > 0$, $r \sim_O r'$ iff $0 \ll r/r' \ll \infty $.
Note also that $0 \sim_O r$ iff $r = 0$. 
If $r_i \sim_O r_1$, where $\{r_i\}_{1 \leq i \leq n}$, is a limited
sequence of non-negative hyperreals, then $\sum_{i=1}^n r_i \sim_O r_1$.

Note that order of magnitude comparisons for individual real numbers
are essentially vacuous.  Similarly, for a function $f$ defined on a
finite interval of $\mathbb{R}$, one cannot usefully assign an order
of polynomial growth to $f$.  However, for functions on intervals in
$\transfer{\mathbb{R}}$ the perspective of nonstandard analysis allows
finer distinctions.

\begin{definition}
Suppose $0 < I_{-} \leq I_{+}$. A $\transfer{\mathbb{R}}$-valued function $f$
is of {\em polynomial growth of order $\lambda$ and scale factor
$\SCALE$ in the interval $[I_{-},I_{+}]$} iff $f(r) \sim_O \SCALE \,
r^\lambda$ for $r \in [I_{-},I_{+}]$.
\end{definition}
The scale factor and the order of growth are not uniquely determined.
However:
\begin{prop}\label{uniqueness-of-growth}
If $f$ has polynomial growth on $[I_{-}, I_{+}]$ with
$I_{-} \preceq_o I_{+}$, then the order of growth is
uniquely determined up to $\cong$.  If there is a hyperreal $\theta$
such that $0 \ll \theta \ll \infty$ and
$I_{-} \leq {(I_{+}/I_{-})}^\theta \leq I_{+}$,
then the scale factor is unique up to $\sim_O$.
\end{prop}
{\sc Proof.} Suppose $f(r) \sim_O \SCALE \,
r^\lambda \sim_O \SCALE' \, r^\mu$ for $r \in [I_{-},I_{+}]$.  Assume without
loss of generality that $\mu > \lambda$.  Thus
\begin{equation}\label{two-scale-comparizon-formula}
\fraco{\SCALE}{\SCALE'} \sim_O r^{\mu - \lambda}
\end{equation}
for all $r \in [I_{-}, I_{+}]$ and so 
$I_{+}^{\mu - \lambda} \sim_O \fraco{\SCALE}{\SCALE'} \sim_O
I_{-}^{\mu - \lambda}$.
Therefore,
\begin{equation}\label{comparison-of-extremes}
{\biggl(\frac{I_{+}}{I_{-}} \biggr)}^{\mu -\lambda} \sim_O 1
\end{equation}
from which follows $(\mu - \lambda) \, \ln (I_{+}/I_{-}) \ll \infty$.  Since
$\ln (I_{+}/I_{-}) \cong \infty$, we conclude that $\mu - \lambda \cong 0$.

If $I_{-} \leq {(I_{+}/I_{-})}^\theta
\leq I_{+}$, then instantiating $r$
in~(\ref{two-scale-comparizon-formula}) with
${(I_{+}/I_{-})}^\theta$, and
using~(\ref{comparison-of-extremes}) and the fact $0 \ll \theta \ll
\infty$, we deduce $\fraco{\SCALE}{\SCALE'} \sim_O
{(I_{+}/I_{-})}^{\theta \, (\mu - \lambda)} \sim_O
1$.~$\blacksquare$

\subsection{Internal Metric Spaces}

Suppose $(\mathfrak{X},d)$ is an internal metric space. For $x,y \in
\mathfrak{X}$, define $x \cong y$ iff $d(x,y)$ is infinitesimal.
``$\cong$'' is an (external) equivalence relation on $\mathfrak{X}$
and $\lhull{\mathfrak{X}}$ is the quotient space
$\mathfrak{X}/\cong$. For $x \in \mathfrak{X},$ let $\hat{x}$ be the
$\cong$-equivalence class of $x$ in $\mathfrak{X}$. We will denote the
canonical quotient map $x \rightarrow \hat{x}$ by
$\varphi_\mathfrak{X}$ and refer to it as the {\em infinitesimal
identification map}. The quotient space $\lhull{\mathfrak{X}}$ becomes
an extended metric space with the extended metric
$\lhull{d}(\hat{x},\hat{y}) = \opr{st}(d(x,y))$.
$\lhull{\mathfrak{X}}$ with the metric $\lhull{d}$ is called the {\em
infinitesimal hull} of $(\mathfrak{X},d)$.
In general $\lhull{\mathfrak{X}}$ is an extended metric space.  
\begin{example}
If $N \cong \infty$, let $J=\{0, h , \ldots, (N-1)h \}$,
where $h=\frac{1}{N}$.  The map $\opr{st}: J \rightarrow [0,1]$ is
equivalent to the infinitesimal identification map, in the following
sense: there is an isometric map $h: [0,1] \rightarrow \lhull{J}$ such
that $h \circ \opr{st} = \varphi_{J}$.
\end{example}

We also define $x \sim y$ to mean $d(x,y) \ll \infty$.  $\sim$ is also
an external equivalence relation. The equivalence classes of
$\mathfrak{X}$ under ``$\sim$'' are the {\em limited components} of
$\mathfrak{X}$.
Clearly the bounded components of $\lhull{\mathfrak{X}}$ are the
images under the map $\varphi_\mathfrak{X}$ of the limited components
of $\mathfrak{X}$.

\begin{example}
Let $J'=\{-N \, h , -(N-1) \, h, , \ldots, 0, \ldots +(N-1) \, h ,N \,
h\}$, where $h = \frac{1}{\sqrt{N}}$.  In this case $\lhull{J'}$ is
the disjoint union of its bounded components each one which is
isometrically isomorphic to $\mathbb{R}$.
\end{example}

It follows immediately from countable saturation, that any separable
metric space $Y$ is isometrically isomorphic to a {\em subset} of a
metric space $\lhull{\mathfrak{X}}$ with $\mathfrak{X}$ hyperfinite.

$\mathfrak{X} \mapsto \lhull{\mathfrak{X}}$ is a covariant functor
from the category of internal metric spaces and $\ESS$-continuous maps
into the category of metric spaces and uniformly continuous mappings and
$\varphi_\mathfrak{X}: \mathfrak{X} \rightarrow \lhull{\mathfrak{X}}$
is a natural transformation of functors.

We have the following relations between balls in $\mathfrak{X}$ and
$\lhull{\mathfrak{X}}$: If $r' < r$ and both are standard, then
\begin{equation}\label{ball-relation-between-beteen-space-and-hull}
\opr{\overline{B}}(x,r') \subseteq \varphi^{-1}_\mathfrak{X}\bigl(\opr{B}(\hat{x},r)\bigr)
\subseteq \opr{B}(x,r).  \end{equation}
Similarly,
\begin{equation}\label{ball-relation-between-beteen-space-and-hull-closed}
\opr{\overline{B}}(x,r') \subseteq
\varphi_\mathfrak{X}^{-1}\bigl(\opr{\overline{B}}(\hat{x},r')\bigr)
\subseteq \opr{B}(x,r).  \end{equation}
For instance~(\ref{ball-relation-between-beteen-space-and-hull-closed})
follows from the implications $d(x,y) \leq r' \implies \opr{st}\,
d(x,y) \leq r' \implies d(x,y) < r$ which are valid so long as $r' <
r$ are both standard.

The following proposition is well-known.  See for
instance~\cite{albeverio-etal}.

\begin{prop}\label{collapsingcomplete} If $(\mathfrak{X},d)$ is an
internal metric space, then each of the bounded components of
$\lhull{\mathfrak{X}}$ is complete. In particular,
$\lhull{\mathfrak{X}}$ is complete. If $C \subseteq \mathfrak{X}$ is
internal, then $\lhull{C}$ (the image of $C$ in
$\lhull{\mathfrak{X}}$) is closed in $\lhull{\mathfrak{X}}$.
\end{prop}
For a metric space $X$ in ${\MDLM}$, define $\iota_X$ so that the
following diagram commutes:
\begin{diagram}
X 	&	&	 \\
\dTo<{\transfer} & \rdTo^{\iota_X} & \\
\transfer X  & \rTo^{\varphi_{\transfer X}}	&  \phull{\transfer X}
\end{diagram}
\begin{prop} \label{embedding-metric} \label{collapsingextends} 
If $X$ is a ${\MDLM}$ metric space, then 
$\iota_X:x \mapsto \varphi_{\transfer X} (\transfer x)$ is an isometric map $X
\rightarrow \phull{\transfer X}$.  
If $X$ is compact, then $\iota_X$ is an isometric isomorphism $X
\rightarrow \phull{\transfer X}$.
\end{prop}

In particular, fix a (non-extended) metric space $(X,d)$.  Since it is
cumbersome to explicitly indicate the $\iota_X$ map $X \rightarrow
\phull{\transfer X}$, we usually view $\iota_X$ is an inclusion map
with the property that $\varphi_{\transfer{X}}(\transfer x) = x$.  Thus $\tilde{X}
= \phull{\transfer X}$ contains $X$ itself and unless $X$ is compact,
$\tilde{X}$ is much larger than $X$. For example, consider
$\mathbb{Z}$ with the metric $d(x,y) =1$ for $x \neq y$.  The ball
$\overline{\opr{B}}(0,1) = \mathbb{Z}$ is not compact since the
sequence $x_n = n$ has no convergent subsequences.  Notice
$\tilde{\mathbb{Z}}$ in this metric has one bounded component which is
(much) larger than $\mathbb{Z}$.  We now determine some conditions
under which the bounded component of $X$ in $\tilde{X}$ is $X$ itself.

\begin{prop} \label{characterization-of-self-bounded-component}
Suppose $(X,d)$ is separable.  A necessary and sufficient condition
that the bounded component of $X$ in $\tilde{X}$ coincide with $X$ is
that every closed ball of $(X,d)$ be compact.
\end{prop}
{\sc Proof.} Observe first that if $a \in X$ is such that
$\opr{\overline{B}}(a,r)$ is compact, then
$\opr{\overline{B}}_{\tilde{X}}(a,r)= \opr{\overline{B}}_X(a,r)$.
Obviously, $\opr{\overline{B}}_{\tilde{X}}(a,r) \supseteq
\opr{\overline{B}}_X(a,r)$.  Let $\varphi$ be the infinitesimal
identification map $\transfer X \rightarrow \tilde{X}$.  If $r < r'$,
by formula~(\ref{ball-relation-between-beteen-space-and-hull-closed})
and well-known nonstandard characterizations of
compactness~(\cite{albeverio-etal}, 2.1.6),
$$\varphi^{-1}\bigl(\opr{\overline{B}}_{\tilde{X}}(a,r)\bigr)
\subseteq \opr{\overline{B}}(\transfer a,r') = \transfer [
\opr{\overline{B}}(a,r')] \subseteq
\varphi^{-1}\bigl(\opr{\overline{B}}_{X}(a,r')\bigr).$$
Since $\varphi$ is surjective, $\opr{\overline{B}}_{\tilde{X}}(a,r)
\subseteq \opr{\overline{B}}_{X}(a,r')$, and as $r'$ is an arbitrary
standard real $> r$, $\opr{\overline{B}}_{\tilde{X}}(a,r) \subseteq
\opr{\overline{B}}_{X}(a,r)$.

To show sufficiency, suppose $a \in X$ and let $x \in \tilde{X}$ be in
the bounded component of $a$. There is a standard $r$ such that
$d_{\tilde{X}}(a, x) \leq r$. By the remarks of the first paragraph,
$\opr{\overline{B}}_{\tilde{X}}(a,r)= \opr{\overline{B}}_X(a,r)$, so
$x \in X$.  Necessity: Suppose some ball $\opr{\overline{B}}_X(a,r)$
with $r \in \mathbb{R}$ is not compact; let $\{y_i\}_{i \in
\mathbb{N}}$ be a sequence of $\opr{\overline{B}}_X(a,r)$ such that
such that $\inf_{i \neq j} d(y_i, y_j) = \rho > 0$.  $y_i = \hat{x_i}$
for some sequence $\{x_i\}_{i \in \mathbb{N}}$ of $\transfer X$.  By
saturation this sequence extends to an internal sequence $x_1, \ldots
, x_N$ and by overspill we can also assume $d(x_i, x_j) \geq \rho/2$
for $1 \leq i,j \leq N$ with $i \neq j$ and $d(\transfer a,x_i) \leq 2
\, r$ for all $i$.  In particular, $d(\hat{x_i}, \hat{x_j}) \geq
\rho/2$ for $1 \leq i,j \leq N$ with $i \neq j$ and $d(a, \hat{x_i})
\leq 2 \, r$ for all $i$. By assumption, the bounded component of $X$
in $\tilde{X}$ is $X$ itself. Thus $\hat{x_i} \in X$ for all $1 \leq i
\leq N$.  Now $N \cong \infty$, so $\{1, \ldots, N\}$ is
uncountable.~$\blacksquare$

\begin{definition}
An internal mapping $f$ from an internal metric space $(\mathfrak{X},
d_{\mathfrak{X}})$ to $(\mathcal{Y}, d_{\mathcal{Y}})$ is
$S$-Lipschitz iff there is a limited hyperreal $K$ such that for all
$x, x' \in \mathfrak{X}$, $d_{\mathcal{Y}}(f(x),f(x')) \leq K \,
d_{\mathfrak{X}}(x,x')$.
\end{definition}
There are related concepts such as $S$-Lipschitz equivalent metrics
whose formulation we leave to the reader.

Finally, we note that the concept of length space is an internal one,
so is applicable to standard metric spaces and internal metric spaces.

\subsection{Borel Structure and Loeb Measure}

We recall two basic concepts in nonstandard measure theory.  The {\em
inner Borel algebra} of an internal set $\mathfrak{X}$ is the smallest
$\sigma$-algebra containing all internal subsets of $\mathfrak{X}$.  A
{\em hypermeasure} on $\mathfrak{X}$ is an internal, nonnegative and
hyperfinitely additive function defined on all internal subsets of
$\mathfrak{X}$.  The conventional name given in the literature is the
less mellifluous ``hyperfinitely additive nonnegative measure''.  The
following result is essentially due to Loeb~\cite{Loeb75} in the case
the hypermeasure is limited.  The uniqueness result in the unlimited
case is due to Henson (Corollary 1 of~\cite{henson}).
\begin{thm} \label{extending} Suppose $\mu$ is a hypermeasure on
$\mathfrak{X}$. Then there is a unique countably additive measure on
the inner Borel sets which extends the mapping $A \mapsto
\opr{st}(\mu(A))$ on the internal sets. 
\end{thm}
The completion of the countably additive measure specified in the
previous result is called the {\em Loeb measure} associated with
$\mu$.   We denote Loeb measure by ${\opr{L}(\mu)}$.

\begin{prop} \label{infinitesial-identification-is-borel} Suppose
$(\mathfrak{X},d)$ is an internal metric space.  Then the
infinitesimal identification map $\varphi_{\mathfrak{X}}: \mathfrak{X}
\rightarrow \lhull{\mathfrak{X}}$ is measurable where $\mathfrak{X}$
has the inner Borel structure and $\lhull{\mathfrak{X}}$ has the ball
Borel structure. \end{prop}
{\sc Proof.}  If $r \in {\mathbb{R}}$ then
$
\hat{y} \in \opr{B}(\hat{x},r) \iff \exists m \in {\mathbb{N}} \
d(x,y) < r-\frac{1}{m} 
\iff y \in \bigcup_{m \in {\mathbb{N}}} \opr{B}(x,r-1/m)
$.
Since each $\opr{B}(x,r-1/m)$ is internal, it follows
$\varphi_\mathfrak{X}^{-1} \opr{B}(\hat{x},r)\}$ is inner
Borel.~$\blacksquare$

\begin{definition}
Suppose $\mathfrak{X}$ is an internal metric space. If $\mu$ is a
hypermeasure on $\mathfrak{X}$, then $\opr{\overline{P}}(\mu)$ is the
measure $\varphi_\ast \opr{L}(\mu)$ defined on the $\sigma$-algebra of
$A \subseteq \lhull{\mathfrak{X}}$ such that $\varphi^{-1}(A)$ is
$\opr{L}(\mu)$ measurable.  $\opr{P}(\mu)$ is the restriction of
$\opr{\overline{P}}(\mu)$ to the ball Borel sets of
$\lhull{\mathfrak{X}}$.
\end{definition}
Note that since any inner Borel set is Loeb measurable,
Propostion~\ref{infinitesial-identification-is-borel} immediately
implies the ball Borel sets of $\lhull{\mathfrak{X}}$ are
$\opr{\overline{P}}(\mu)$ measurable.

Since $\opr{L}(\mu)$ is complete, $\opr{\overline{P}}(\mu)$ is also
complete.  Moreover, $\opr{\overline{P}}(\mu)$ is not too far off from
being the completion of $\opr{P}(\mu)$.

\begin{prop}\label{completion-of-reduction}
Suppose $\mu$ is a hypermeasure on $\mathfrak{X}$ such that
$\opr{P}(\mu)(K) < \infty$ for every compact $K \subseteq
\lhull{\mathfrak{X}}$.  If $A \subseteq \lhull{\mathfrak{X}}$ is a
subset of a $\sigma$-compact set and $A$ is
$\opr{\overline{P}}(\mu)$-measurable then $A$ is in the completion of
$\opr{P}(\mu)$.  Conversely, any set in the completion of
$\opr{P}(\mu)$ is $\opr{\overline{P}}(\mu)$ measurable.
\end{prop}
{\sc Proof.}  Suppose $A$ is $\opr{\overline{P}}(\mu)$-measurable.  By
assumption, $A \subseteq \bigcup_n K_n$ with $K_n$ compact.  Since the
domain of the completion of a measure is a $\sigma$-algebra, it
suffices to show each $A \cap K_n$ is in the completion of
$\opr{P}(\mu)$.  Thus without loss of generality, we can assume $A
\subseteq K$ for some compact $K$.  In particular, 
$$\opr{L}(\mu) ( \varphi_{\mathfrak{X}}^{-1}(A) = \overline{\opr{P}}(\mu) (A) \leq
\opr{P}(\mu) (K) < \infty.$$

Let $\epsilon > 0$ with $\epsilon \in \mathbb{R}$ be arbitrary and $C
\subseteq \varphi_{\mathfrak{X}}^{-1}(A)$ internal such that
$\opr{L}(\mu) ( \varphi_{\mathfrak{X}}^{-1}(A) \setminus C ) <
\epsilon$.  $\varphi_{\mathfrak{X}}(C)$ is closed and
$\varphi_{\mathfrak{X}}(C) \subseteq A \subseteq K$ so
$\varphi_{\mathfrak{X}}(C)$ is compact and hence ball Borel.
Moreover,
\begin{eqnarray*}
\opr{\overline{P}}(\mu)(A \setminus \varphi_{\mathfrak{X}}(C)) & = &
\opr{L}(\mu)\bigl(\varphi_{\mathfrak{X}}^{-1}(A \setminus \varphi_{\mathfrak{X}}(C))\bigr) \\
& = & \opr{L}(\mu)\bigl(\varphi_{\mathfrak{X}}^{-1}(A) \setminus
\varphi_{\mathfrak{X}}^{-1}(\varphi_{\mathfrak{X}}(C))\bigr) \\
& \leq &  \opr{L}(\mu)(\varphi_{\mathfrak{X}}^{-1}(A) \setminus C) \leq \epsilon
\end{eqnarray*}
Since $\epsilon > 0$ is arbitrary, it follows there is a
$\sigma$-compact (and therefore ball Borel set) $B \subseteq A$ such
that $\opr{\overline{P}}(\mu)(A \setminus B) = 0$. Applying this to $K
\setminus A$, we conclude there is also a Borel set $B'$ such that $A
\subseteq B' \subseteq K$ and $\opr{\overline{P}}(\mu)(B' \setminus A)
= 0$. Thus $A$ is nested between ball Borel sets $B$ and $B'$ with
$\opr{\overline{P}}(\mu)(B' \setminus B) = 0$, and so $A$ is in the
completion of $\opr{P}(\mu)$. The converse follows from the
definitions.~$\blacksquare$

\section{Compactness and Measure}

We begin by giving some general compactness properties which follow
from the existence of a Borel measure.  All statements of this section
are standard.  If $\nu$ is a countably additive measure on a
measurable space $(X. \mathcal{A})$, a set $A \subseteq X$ measurable
or not is $\nu$-finite iff there is a $B \in \mathcal{A}$ such that $A
\subseteq B$ and $\nu(B) < \infty$.  Any subset of a $\nu$-finite set
is $\nu$-finite.  In the context of the following result
note~\cite{christensen}.

\begin{prop}\label{characterization-of-compactness}
Let $(X,d)$ be a complete extended metric space and $\nu$ a countably
additive ball Borel measure on $X$ with the following property: There
is an $R > 0$ such that
$$\nu\bigl(\opr{B}(x,r)\bigr) < \infty \quad\mbox{
and
 } \quad 0 < \inf_{x,x'}
\frac{\nu\bigl(\opr{B}(x,r)\bigr)}{\nu\bigl(\opr{B}(x',r)\bigr)}$$
for all $0 < r <R$.  Then $K \subseteq X$ is relatively compact iff
there is a $\delta > 0$ such that $\{x: d(x,K) \leq \delta\}$ is
$\nu$-finite.
\end{prop}
{\sc Proof.} Suppose $K \subseteq X$ is relatively compact. Let
$\epsilon \in \mathbb{R}$ be such that $0 < \epsilon < R$.  By
relative compactness of $K$, there is an $n \in \mathbb{N}$ and a
sequence ${x}_1, \ldots, {x}_n \in K$ such that
$\overline{K} \subseteq \bigcup_{i=1}^n \opr{B}({x}_i,\epsilon) = W$.
$W$ is open, so
$
\delta = \inf \{d(x,y): y \in {X} \setminus W, x \in \overline{K}\} >
0$.
Let $V=\{x: d(x,K) \leq \delta/2\}$. Clearly $V \subseteq W$. Now
$\bigcup_{i=1}^n \opr{B}(x_i,\epsilon)$ is ball Borel and
$\nu$-finite. In particular $V$ is $\nu$-finite, proving the claim.
%
%

Conversely, suppose $\delta>0$ is such that $V = \{x: d(x,K) \leq
\delta\}$ is
$\nu$-finite.
We may assume $\delta < R$.  We prove by contradiction that $K$ is
precompact.  By possibly choosing a smaller positive $\delta$, we may
assume
there is an infinite sequence $\{x_i\}_{i \in \mathbb{N}}$ in $K$ such
that $x_i \not\in \opr{B}(x_j,\delta)$ for $j<i$.  Thus
$\{\opr{B}(x_i,\delta/2)\}_{i \in \mathbb{N}}$ is a sequence of
pairwise disjoint balls.  Furthermore, $\opr{B}(x_i, \delta/2)
\subseteq V$, since $x_i \in K$. Thus for any $N$
$$\begin{aligned}
\nu \bigl(\bigcup_{i=1}^\infty \opr{B}(x_i, \delta/2)\bigr) = &
\sum_{i=1}^\infty \nu(\opr{B}(x_i, \delta/2)) \\
\geq & N \ \nu(\opr{B}(x_1, \delta/2)) \ \inf_{1 \leq i \leq \infty}
\frac{\nu(\opr{B}(x_i, \delta/2))}{\nu(\opr{B}(x_1, \delta/2))}
\end{aligned}
$$
so that $\nu \bigl(\bigcup_{i=1}^\infty \opr{B}(x_i, \delta/2)\bigr) =
+\infty$.  Since $V \supseteq \bigcup_{i=1}^\infty \opr{B}(x_i,
\delta/2)$, $V$ cannot be $\nu$-finite.  It follows that $K$ is
precompact.  Thus by completeness of $(X,d)$ the closure of $K$ is
compact.~$\blacksquare$

\begin{cor}\label{characterization-of-localcompactness} 
Let $(X,d)$ be a complete extended metric space and $\nu$ a ball Borel
measure on $X$ satisfying the conditions of
Proposition~\ref{characterization-of-compactness}. Then $X$ is locally
compact.  In fact, all closed balls of radius $< R$ are
compact.\end{cor}
{\sc Proof.} Suppose $x \in X$ and $r$ is a real with $0< r <
R$. Consider the set
$V=\{y \in X: d(y,\opr{B}(x,r)) \leq \delta\} \subseteq \opr{B}(x,r + \delta)$.
If $r + \delta < R$, then $V$ is $\nu$-finite. Thus by
Proposition~\ref{characterization-of-compactness}, $\opr{B}(x,r)$ is
relatively compact.~$\blacksquare$

\section{Near Homogeneity}

\begin{definition} \label{near-homogeneity-definition}
A hypermeasure $\mu$ on an internal metric space $(\mathfrak{X},d)$ is
{\em nearly homogeneous} iff there is an $R_\mu \gg 0$ such
that for all open balls $V$, $V'$ of radius $s, s'$ respectively with
$0 \ll s ,s' \leq R_\mu$, $\mu(V) \sim_O \mu(V')$ and
$\mu(V) \neq 0$.  $R_\mu$ is a {\em radius of homogeneity}
of $\mu$.
\end{definition}
In this definition, note the possibility that
$R_\mu$ is an unlimited hyperreal.  

Denote the common $\sim_O$-equivalence class of the hyperreals
$\mu(\opr{B}(x,s))$, for $0 \ll s \leq R_\mu$ by
$M_\mu$.  Note that this is an external set.  Elements of
$M_\mu$ are normalization constants of $\mu$.  By abuse of
notation, we write $r \sim_O M_\mu$ instead of $r \in
M_\mu$.  Similarly, we write $r \preceq_O M_\mu$
iff for any $s \in M_\mu$, $r \preceq_O s$.  If
$M_\mu \sim_O 1$, we say $\mu$ is {\em normalized}.

Given a nearly homogeneous hypermeasure $\mu$ on an internal metric
space $(\mathfrak{X},d)$, we associate to it a class of normalized
hypermeasures on $\mathfrak{X}$ as follows: For any $M \sim_O
M_\mu$, let
$$\mu_M(A) = \frac{1}{M} \mu(A).$$
The hypermeasures $\mu_M$ are all essentially equivalent. For
instance, the measures $\mu_M$ are constant multiples $r$, $0 \ll r
\ll \infty$, of each other. It follows that the countably additive
measures $\opr{L}(\mu_M)$ have identical null sets and sets of finite
measure. 

We will call an internal metric space $(\mathfrak{X},d)$ nearly
homogeneous if the hypermeasure $\mu_{\opr{card}}: A \mapsto
\opr{card} A$ is nearly homogeneous, where we make the convention that
$\mu_{\opr{card}}(A) = + \infty \in \transfer \overline{\mathbb{R}}$
in case $A \subseteq \mathfrak{X}$ is not hyperfinite. Henceforth, if
$(\mathfrak{X},d)$ nearly homogeneous, we will use $\mu_{\opr{card}}$
to denote a normalized version of $A \mapsto \opr{card}(A)$.

In particular:

\begin{prop}\label{near-homogeneous-measures-of-finite-balls}
Suppose $\mu$ is a normalized nearly homogeneous hypermeasure on
$(\mathfrak{X},d)$.  Then the measure $\opr{P}(\mu)$ on
$\lhull{\mathfrak{X}}$ has the property that
\begin{equation} \label{compactness-measure-condition-1}
0 < \opr{P}(\mu)\bigl(\opr{B}(x,r)\bigr) < \infty
\end{equation}
and
\begin{equation} \label{compactness-measure-condition-2} 
0 < \inf_{x,x'}
\frac{\opr{P}(\mu)\bigl(\opr{B}(x,r)\bigr)}{\opr{P}(\mu)\bigl(\opr{B}(x',r)\bigr)}
\end{equation}
for all $r \in \mathbb{R}$ with $0 < r \leq R_\mu$.
\end{prop}
{\sc Proof.}  Let $\varphi: \mathfrak{X} \rightarrow
\lhull{\mathfrak{X}}$ be the infinitesimal identification map. Thus
$$\opr{B}(x,\frac{r}{2}) \subseteq
\varphi^{-1}\bigl(\opr{B}(\hat{x},r)\bigr) \subseteq \opr{B}(x, r).$$
Since $\mu$ is normalized, $\mu \bigl(\opr{B}(x,\frac{r}{2})\bigr)
\sim_O \mu \bigl(\opr{B}(x, r)\bigr) \sim_O 1$, so
$$0 \ll \opr{st} \mu \bigl(\opr{B}(x,\frac{r}{2})\bigr) \leq
\opr{P}(\mu)\bigl(\opr{B}(\hat{x},r)\bigr) \leq \opr{st} \mu
\bigl(\opr{B}(x, r)\bigr) \ll \infty.$$
This proves~\ref{compactness-measure-condition-1}. To
prove~\ref{compactness-measure-condition-2}, note that
$$\frac{\opr{P}(\mu)\bigl(\opr{B}(\hat{x},r)\bigr)}
{\opr{P}(\mu)\bigl(\opr{B}(\hat{x}',r)\bigr)} \geq \opr{st}
\biggl(\frac{\mu \bigl(\opr{B}(x,\frac{r}{2})\bigr)}{\mu
\bigl(\opr{B}(x', r)\bigr)}\biggr).$$
However, for each $x,x' \in \mathfrak{X}$ and each standard $r$ for which $0 < r
\leq R_\mu$,
$$\frac{\mu \bigl(\opr{B}(x,\frac{r}{2})\bigr)}{\mu \bigl(\opr{B}(x',
r)\bigr)} \gg 0.$$
For a fixed value of $r$ the above expression is an internal function
of $(x,x')$. Therefore for a fixed standard value of $r$, its infimum
over all pairs $(x,x')$ is $\gg 0$
and~(\ref{compactness-measure-condition-2}) follows.~$\blacksquare$

In particular:

\begin{prop}\label{nearly-homogeneous-implies-locally-compact}
Suppose $\mu$ is a normalized nearly homogeneous hypermeasure on
$(\mathfrak{X},d)$. All closed balls of radius $<
R_\mu$ in $\lhull{\mathfrak{X}}$ are
compact.
\end{prop}

\section{Polynomial Growth} \label{polynomial-growth-section}

\begin{definition}
Let $(\mathfrak{X},d)$ be an internal metric space, $I_{-} \leq I_{+}$
nonnegative hyperreals. $(\mathfrak{X},d)$ is of {\em polynomial
growth in the interval $[I_{-},I_{+}]$} iff there is a hyperreal $\SCALE$ and
a hyperreal $0 \ll \lambda \ll \infty$ such that for $I_{-} \leq r \leq
I_{+}$
\begin{equation} \label{polynomial-growth-equation}
\opr{card} \, \opr{B}(x,r) \sim_O \SCALE \, r^\lambda.  
\end{equation}
\end{definition}
The hyperreal $\SCALE$ is the {\em scale factor} and the limited
hyperreal $\lambda$ is the {\em order of growth} of $\mathfrak{X}$ in
$[I_{-},I_{+}]$.  Alternatively, $(\mathfrak{X},d)$ is of polynomial growth of
order $\lambda$ in $[I_{-},I_{+}]$ with scale factor $\SCALE$ iff there
are constants $c \gg 0$ and $C \ll \infty$ such that
\begin{equation} \label{altrenate-polynomial-growth-equation}
c \, \SCALE \, r^\lambda \leq \opr{card} \, \opr{B}(x,r) \leq C \,
\SCALE \, r^\lambda.
\end{equation}
for $r \in [I_{-}, I_{+}]$.  We will refer to the constants $c$, $C$ as the
lower and upper bounds of polynomial growth of $(\mathfrak{X},d)$ in $[I_{-}, I_{+}]$.

Unless otherwise stated, we will implicitly assume $I_{-} \preceq_o I_{+}$.
The following immediately follows from
Proposition~\ref{uniqueness-of-growth}.
\begin{prop}
If $I_{-} \preceq_o I_{+}$, then the order of growth is uniquely
determined up to $\cong$.  If there is a hyperreal $\theta$ such that
$0 \ll \theta \ll \infty$ and $I_{-} \leq
{(I_{+}/I_{-})}^\theta \leq I_{+}$, the scale factor is unique up to
$\sim_O$.
\end{prop}

For each nonnegative $\Rescale \in \transfer\mathbb{R}$, let $d_\Rescale(x,y) =
\frac{1}{\Rescale}d(x,y)$. We use the notation $\opr{B}_\Rescale$ to denote open
balls relative to $d_\Rescale$. Clearly, $d_\Rescale(x,y) < r$ iff $d(x,y) < \Rescale r$ so
that
$\opr{B}_\Rescale(x,r)=\opr{B}(x, \Rescale r)$.

Consider $\mathfrak{X}$ as a hypermeasure space with counting measure $\mu(V) =
\opr{card}(V)$.

\begin{prop}
Suppose $I_{-} \preceq_o I_{+}$ and $\Rescale > 0$ is such that
$I_{+}/\Rescale \gg 0$ and $\Rescale/I_{-} \cong +\infty$.  If
$(\mathfrak{X},d)$ is of polynomial growth in $[I_{-}, I_{+}]$ then
counting measure $\mu$ is nearly homogeneous on
$(\mathfrak{X},d_\Rescale)$ with a radius of homogeneity $R
=I_{+}/\Rescale $ and normalizing constant $M_\mu = \SCALE \,
\Rescale^\lambda$.

In particular, closed balls of radius $<R$ are compact in
$\lhull{(\mathfrak{X},d_\Rescale)}$ and $\lhull{\mathfrak{X}}$ is locally
compact.
\end{prop}
{\sc Proof.} If $0 \ll r,r' \leq R$, then $I_{-} \leq r \Rescale,r' \Rescale \leq R
\, \Rescale = I_{+}$, so
$$\opr{card} \, \opr{B}_\Rescale(x,r)  = \opr{card} \, \opr{B}(x, \Rescale r) \sim_O
\SCALE \, (\Rescale r)^\lambda \sim_O \SCALE \, (\Rescale r')^\lambda
 = \opr{card} \, \opr{B}_\Rescale(x,r').$$
By Proposition~\ref{characterization-of-localcompactness}, local
compactness of $\lhull{\mathfrak{X}}$ follows.~$\blacksquare$

If $(\mathfrak{X},d)$ is an internal metric space of polynomial growth
$\lambda$ in $[I_{-}, I_{+}]$, then $\mu(A) =
\opr{card} A/\SCALE \, \Rescale^\lambda$ is a normalized homogeneous
hypermeasure on $(\mathfrak{X}, d_\Rescale)$.  $\lhull{(\mathfrak{X}, d_\Rescale)}$
is locally compact and is equipped with a countable additive ball
Borel measure $\opr{P}(\mu)$ which is finite on compacts.

\begin{thm} \label{measure-characterization-polynomial-growth} Suppose
$(\mathfrak{X},d)$ is of polynomial growth $\lambda$ in
$[I_{-}, I_{+}]$ and $\Rescale>0$ is such that $R =
I_{+}/\Rescale \gg 0$ and $\Rescale/I_{-} \cong +\infty$.
The countably additive measure $\opr{P}(\mu)$ on
$\lhull{(\mathfrak{X}, d_\Rescale)}$ has the property that there are
constants $k,K \in \mathbb{R}$ such that
$$k \, r^{\opr{st}(\lambda)} \leq \opr{P}(\mu)\bigl(\opr{B}(x,r)\bigr)
\leq K \, r^{\opr{st}(\lambda)}$$
for every $x \in \lhull{\mathfrak{X}}$ and $r \in \mathbb{R}$ such that $0
< r \leq R$.  
Note that $\opr{B}(x,r)$ refers to the ball in the metric $\lhull{d_\Rescale}$.
\end{thm}
{\sc Proof.} By definition, there are standard constants $0 < c \leq C
< \infty$ such that for $0 \ll r \leq R$, $x \in \mathfrak{X}$
$$
c \leq \frac{\opr{card} \, \opr{B}_\Rescale(x,r)}{\SCALE \, (r \Rescale)^\lambda} \leq C
$$
Now use formula~(\ref{ball-relation-between-beteen-space-and-hull}) on
$(\mathfrak{X},d_\Rescale)$.  Thus if $r$ is a nonnegative standard real such
that $r\leq R$,
$$
c \, \biggl(\frac{r}{2}\biggr)^{\opr{st}(\lambda)}  \leq \opr{st} \,
\mu\bigl(\opr{B}_\Rescale(x,\frac{r}{2})\bigr) \leq
\opr{P}(\mu)\bigl(\opr{B}(\hat{x},r)\bigr) 
\leq \opr{st}
\mu\bigl(\opr{B}_\Rescale(x, r)\bigr) 
\leq C r^{\opr{st}(\lambda)}.
$$

\begin{example}
Consider a hyperfinite rooted tree $T$ with leaf nodes
$\opr{leaf}({T})$. Any $x \in \opr{leaf}({T})$ in $T$ is uniquely
identified by the sequence of branches
$\opr{path}(x) = \{\opr{b}_i(x)\}_{1 \leq i \leq \opr{depth}(x)}$
required to reach $x$ from the root node of $T$.
If $m$ is a standard integer, define a metric $d_m$ on
$\opr{leaf}({T})$ as follows: If $x \neq y$,
$d_m(x,y)=m^{-\opr{v}(x,y)}$ where $\opr{v}(x,y)$ is the first index
at which $\opr{path}(x),\opr{path}(y)$ differ, otherwise,
$d_m(x,y)=0$.
It is well-known $d_m$ is an ultrametric.

\begin{prop}
If ${T}$ is a tree of uniform depth $N$ all of whose non-leaf nodes
have standard branching degree $n$, then for $r \in ]m^{-(N+1)}, 1]$.
$$n^{N - 1 } \, r^{\ln n /\ln m} \leq \opr{card} \, \opr{B}(x,r) \leq n^N
\, r^{\ln n /\ln m}.$$
In particular, $(\opr{leaf}({T}), d_m)$ is of polynomial growth of
order $\ln n/\ln m$ in the interval $]m^{-(N+1)}, 1]$.
\end{prop}
{\sc Proof.} 
%
%
Note that the condition $ r \in ]m^{-(N+1)}, 1]$ is equivalent to
\begin{equation} \label{translated-range-condition}
0 \leq  - \ln r / \ln m < N+1.
\end{equation}
If $x \in \opr{leaf}({T})$, then $\opr{B}(x,r)$ consists of $x$ and
those $y \in \opr{leaf}({T})$ with $x \neq y$ and $\opr{v}(x,y) > -
\ln r / \ln m$;
$\opr{B}(x,r)$ is thus the set of those $y$ for which $\opr{b}_i(x) =
\opr{b}_i(y)$ for $1 \leq i \leq  - \ln r / \ln m$.  Counting the
number of all paths (from root to leaf node) which agree with
$\opr{path}(x)$ up to index $ \lfloor - \ln r / \ln m \rfloor$
using~(\ref{translated-range-condition}), it follows that
$$\opr{card} \, \opr{B}(x,r) = n^{N - (\lfloor - \ln r / \ln m \rfloor)}.$$
%
The result now follows from the estimate
$$ n^{N+\ln r/\ln m -1} \leq \opr{card} \, \opr{B}(x,r) \leq n^{N+\ln
r/\ln m},$$
and the fact that for any $r$,
$n^{\ln r/\ln m} = r^{\ln n/\ln m}$.~$\blacksquare$

Note that the Cantor middle thirds set has a Lipschitz equivalent leaf
metric with $m=3$ and $n=2$. We also alert the reader to the fact that
the spaces in this example are not length spaces, since they are
ultrametric spaces and therefore totally disconnected.
\end{example}

\section{Ahlfors Regularity}

All results in this section are standard.  $\mathcal{H}^\lambda$
denotes $\lambda$-dimensional Hausdorff measure, which is a countably
additive measure on the Borel sets, and in particular on the ball
Borel sets.  For definition and background on Hausdorff measures,
see~\cite{mattila},~\cite{rogers}.  These references do not explicitly
treat Hausdorff measures in extended metric spaces, but there is no
difficulty extending the definitions and basic properties to this
context.  Nevertheless, we need to add the following caveat in case
$E$ is a subset of an extended metric space which does not have a
countable $\delta$ cover~(\S 1.2 of~\cite{falconer}).  In this case, we make
the convention that $\mathcal{H}^\lambda_\delta(E) = +\infty$ and
accordingly, $\mathcal{H}^\lambda(E) = +\infty$.

Note that
if $E \subseteq X$ is nonseparable, then $\mathcal{H}^\lambda(E) =
+\infty$ for any $\lambda$.
For if $\mathcal{H}^\lambda < \infty$ then for every $\delta > 0$
there is a countable set $D_\delta\subseteq X$ with the property that
$d(x,D_\delta) \leq \delta$ for every $x \in E$. Thus $D=\bigcup_n
D_{1/n}$ is countable and $d(x, D) = 0$ for all $x \in E$.  In
particular, an uncountable set which intersects each bounded component
in exactly one point has $\lambda$-dimensional Hausdorff measure
$+\infty$ for any $\lambda$. This oddity clearly suggests that we
stick to separable subspaces whenever possible.

It is no great surprise that
Theorem~\ref{measure-characterization-polynomial-growth} determines
the Hausdorff dimension of the extended metric space
$\lhull{(\mathfrak{X},d_\Rescale)}$.

\begin{prop}\label{ahlfors-condition}
Suppose $(X,d)$ is a complete extended metric space and $\nu$ a ball
Borel measure on $X$.  If there are real numbers $R> 0$ and $\lambda >
0$ such that for $0 < r < R$
\begin{equation}
\label{ahlfors-regular-measure-condition}
k \, r^\lambda \leq \nu \bigl(\opr{B}(x,r)\bigr) \leq K \,
r^\lambda
\end{equation}
then all closed balls in $X$ of radius $<R$ are compact. Moreover,
\begin{equation} \label{hausdorff-measure-dominates}
\nu(E) \leq K \ \mathcal{H}^\lambda(E)
\end{equation}
for any ball Borel set $E$.  

For the inequality in the other direction, there is a constant $c > 0$
such that for every Borel set $E$ contained in a separable subset of
$(X, d)$ with $\nu(E) < \infty$, then $\mathcal{H}^\lambda(E) <
\infty$ and
\begin{equation}\label{counting-upper-bound-for-hausdorff-measure}
c \, \mathcal{H}^\lambda (E) \leq \nu(E)
\end{equation}
\end{prop}
{\sc Proof.}  To show compactness of balls of radius $< R$, apply
Corollary~\ref{characterization-of-compactness} to the complete
metric space $(X,d)$ and the ball Borel masure $\nu$.  We verify that
the assumptions on $\nu$ hold:
$\nu\bigl(\opr{B}(x,r)\bigr) \leq K \, r^\lambda < \infty$ and for all
$x,x' \in X$ and all $r < R$,
$$\frac{\nu\bigl(\opr{B}(x,r)\bigr)}{\nu\bigl(\opr{B}(x',r)\bigr)}
\geq k/K > 0.$$

Let $E$ be an arbitrary set, and $\{A_i\}_{i \in \mathbb{N}}$ any
sequence which covers $E$. Then choosing $x_i \in A_i$,
$$\nu(E) \leq \sum_i \nu(A_i) \leq \sum_i \nu\bigl(\opr{B}(x_i,
\opr{diam}A_i)\bigr) \leq K \sum_i {\opr{diam}A_i}^\lambda.$$
Therefore, 
$\nu(E) \leq K \, \mathcal{H}^\lambda_\delta (E)$
for any positive $\delta$ and letting $\delta \rightarrow 0$, it
follows that formula~(\ref{hausdorff-measure-dominates}) holds for any
ball Borel set $E$.

In order to show the inequality in the other direction,
we need a polynomial upper bound on the sizes of coverings by balls.
This is proved using well-known Vitali covering arguments.  See for
example~\cite{falconer}, Theorem 1.10.  For completeness (and since we
have been unable to locate a reference with results in the precise
form we need here,) we provide the details in the following
subsection.

\subsection{Vitali Covering Arguments} We begin with a lemma in which
the cardinality bound on coverings is recast, by brute force, into a
slightly more suggestive format:

\begin{lemma} \label{brute-force-reformulation} 
Let $(X,d)$ be a metric space, $a \in X$, $r>0$ and $\lambda \in ]0,
\infty[$. Suppose there is a $\sigma_0 > 0$ such that for all $\sigma
\leq \sigma_0$, there is an $Y_\sigma \subseteq X$ satisfying
$d(z,Y_\sigma) < \sigma$ for all $z \in X$ and
\begin{equation}\label{exact-dimension-criterion}
C= \sup_{\sigma \in ]0, \sigma_0]} \frac{\opr{card}\bigl(\opr{B}(a,
r)\cap Y_\sigma \bigr)}{{\bigl(\frac{r}{\sigma}\bigr)}^\lambda} <
\infty.
\end{equation}
Then the $\lambda$-dimensional Hausdorff measure of $\opr{B}(a,
\fraco{r}{2})$ is at most $2^\lambda \, C \, r^\lambda$.
\end{lemma}
{\sc Proof.}  Suppose $\sigma \in ]0, \min(\sigma_0, r/2)]$. Every $y
\in \opr{B}(a, r/2)$ is at distance $< \sigma$ from some $z_y \in
Y_\sigma$; moreover by the triangle inequality and the fact $\sigma
\leq r/2$, $z_y \in \opr{B}(a, r)$. By
formula~(\ref{exact-dimension-criterion}) there are at most $C \,
(\fraco{r}{\sigma})^\lambda$ distinct $z_y$.  In particular, the ball
$\opr{B}(a, r/2)$
can be covered by the family $\{\opr{B}(z_y, \sigma)\}_y$ of
cardinality $\leq C (\fraco{r}{\sigma})^\lambda$.  Each such ball has
diameter $ \leq 2 \sigma$.
Thus applying the definition of Hausdorff measure,
$$\mathcal{H}^\lambda_\sigma \biggl(\opr{B}(a, \frac{r}{2})\biggr)
\leq 2^\lambda \, \sigma^\lambda \, C {\biggl(\frac{r}{\sigma}\biggr)}^\lambda
 = 2^\lambda \, r^\lambda \, C.$$
This inequality is valid for $0 < \sigma \leq \min(\sigma_0, r/2)$.
Thus the $\lambda$-dimensional Hausdorff measure of $\opr{B}(a, r/2
)$ is $\leq 2^\lambda \, r^\lambda \, C$ as claimed.~$\blacksquare$

It remains to find a suitable supply of such sets $Y_\sigma$.

\begin{definition} \label{delta-separated-peoperty}
An extended metric space $(X,d)$ is $\sigma$-separated iff $d(x,y)
\geq \sigma$ for all $x,y \in X$ with $x \neq y$.
\end{definition}
For any metric space $(X,d)$, by Zorn's lemma there is a set $Y
\subseteq X$ which is maximal with respect to the property of being
$\sigma$-separated. If $Y \subseteq X$ is maximal $\sigma$-separated,
then for any $x \in X$, $d(x,Y) < \sigma$; otherwise, if $d(x,Y)
\geq \sigma$, then $\{x\} \cup Y$ is also $\sigma$-separated.

\begin{prop} \label{polynomial-growth-implies-robustness}
Suppose $(X,d)$ is a metric space, $\nu$ a ball Borel measure on
$X$.  Let $R> 0$ and $\lambda > 0$ be such that
$\nu\bigl(\opr{B}(a, R)\bigr) < \infty$
and 
$\nu(\opr{B}(x, r)) \geq k \, r^\lambda$
for $x \in \opr{B}(a, R/2)$ and every $r \leq R/2$.  Then if $r \leq
R/2$, the $\lambda$-dimensional Hausdorff measure of $\opr{B}(a,
\fraco{r}{2})$ is at most $4^\lambda \, k^{-1} \, \nu \bigl(\opr{B}(a,
2 \, r ) \bigr)$.
\end{prop}
{\sc Proof.}  Let $r \leq R/2$.  Suppose $0< \sigma \leq r$ and
$Y_\sigma$ is maximal $\sigma$-separated.  We will show
$$D = \frac{\opr{card}\bigl(\opr{B}(a, r)\cap Y_\sigma
\bigr)}{{\bigl(\frac{r}{\sigma}\bigr)}^\lambda} \leq 2^\lambda  \,
k^{-1} \, r^{-\lambda} \, \nu \bigl(\opr{B}(a, 2 \, r ) \bigr) < \infty,$$
independently of $\sigma$.  The result will then follow by
Lemma~\ref{brute-force-reformulation}.  Since $Y_\sigma$ is
$\sigma$-separated, the balls $\opr{B}(z, \sigma/2)$ for $z \in Y$
are pairwise disjoint. Moreover,
$\opr{B}(a, 2 \, r) \supseteq \bigcup_{z \in \opr{B}(a, r)
\cap Y_\sigma} \opr{B}(z, \fraco{\sigma}{2})$.
Thus
$$
\nu \bigl(\opr{B}(a, 2 \, r ) \bigr)  \geq \sum_{z \in \opr{B}(a, r) \cap
Y} \nu\bigl(\opr{B}(z, \frac{\sigma}{2})\bigr)
\geq D \, \biggl(\frac{r}{\sigma}\biggr)^\lambda \ k \,
\biggl(\frac{\sigma}{2}\biggr)^\lambda 
 = 2^{-\lambda} \, D \, k \, r^\lambda
$$
%
%
The estimate on $D$ follows.~$\blacksquare$

Combining Lemmas~\ref{brute-force-reformulation} and
\ref{polynomial-growth-implies-robustness} we very nearly have
$\mathcal{H}^\lambda$ is absolutely continuous with respect to $\nu$
under the assumptions of Proposition~\ref{ahlfors-condition}.  This is
indeed true, as is shown by a Vitali covering argument.  See for
instance~\cite{falconer}.  We state this as a lemma:

\begin{lemma}
Let $(X,d)$ be an extended metric space and $0 < \lambda <
\infty$. Given an open set $U$ in $(X,d)$ and $\rho > 0$, define
inductively a sequence of numbers $\{R_i\}_i$ and closed balls
$V_i=\opr{\overline{B}}(x_i,r_i)$ with the following properties:
\begin{enumerate}
\item $R_{m+1} = \sup \{r \leq \rho: \exists x \in X \
\opr{\overline{B}}(x,r) \subseteq U \setminus \bigcup_{i=1}^m V_i\}$,
where it is understood that if the set of real numbers in braces is
empty then the sequences $\{R_i\}_i$ and $\{V_i\}_i$ terminate at
$m$. Note that if $R_{m+1}$ is defined, then it automatically follows
that $R_{m+1}> 0$.
\item If $R_{m+1} > 0$, let $x_{m+1}$, $r_{m+1}$ be such that $R_{m+1}
\geq r_{m+1} > R_{m+1} /2 $ and $\opr{\overline{B}}(x_{m+1},r_{m+1})
\subseteq U \setminus \bigl(\bigcup_{i=1}^m V_i \bigr)$.
\end{enumerate}
Then either $\sum_i r_i^\lambda = +\infty$ or $\mathcal{H}^\lambda(U -
\bigcup_i V_i) = 0$.
\end{lemma}
{\sc Proof.}  If the sequence $\{V_i\}_i$ is finite of length $m$,
then $U = \bigcup_{i=1}^m V_i$. Suppose $\sum_i r_i^\lambda < \infty$.  We
show that for every $k \in \mathbb{N}$,
\begin{equation} \label{union-decomposition-formula}
U = \bigcup_{i=1}^k V_i \cup \bigcup_{i=k+1}^\infty
\opr{\overline{B}}(x_i, 3 \, r_i).
\end{equation}
To see this, suppose $x \in U \setminus \bigcup_{i=1}^k V_i$.  $U$ is
open, so there is an $r>0$ such that $\opr{\overline{B}}(x,r)
\subseteq U \setminus \bigcup_{i=1}^k V_i$. Since $r_i \rightarrow 0$,
there is an $m$ such that such that $r_m < r/2$.
$\opr{\overline{B}}(x,r)$ must intersect at least one of $V_{k+1},
\ldots, V_{m-1}$.  Otherwise $\opr{\overline{B}}(x,r) \subseteq U
\setminus \bigcup_{i=1}^{m-1} V_i$ so that $R_m \geq r$ and hence $r_m
> R_m/2 \geq r/2$ which contradicts $r_m < r/2$.  Let $\ell$ be the
first index $\geq k+1$ such that $\opr{\overline{B}}(x,r)$ intersects
$V_\ell$.  Thus $R_\ell \geq r$ and $r_\ell > R_\ell/2 \geq r/2$.  By
the triangle inequality
$$d(x_\ell, x) \leq r_\ell + r \leq 3 \, r_\ell,$$
so $x \in \bigcup_{i=k+1}^\infty \opr{\overline{B}}(x_i, 3 \, r_i)$
from which~(\ref{union-decomposition-formula}) follows.  If $\delta > 0$
and $k$ in~(\ref{union-decomposition-formula}) is such that
$\opr{diam}\opr{\overline{B}}(x_i, 3 \, r_i) = 6 \, r_i <
\delta$,  then,
$$\mathcal{H}^\lambda_\delta(U \setminus \bigcup_{i=1}^\infty V_i)
\leq 6^\lambda \, \sum_{i=k+1}^\infty r_i^\lambda \longrightarrow 0.$$
Thus $\mathcal{H}^\lambda_\delta(U \setminus \bigcup_{i=1}^\infty V_i)
= 0$ and so taking the limit as $\delta \longrightarrow 0$, the result
follows.~$\blacksquare$

\paragraph{Completion of Proof of
Proposition~\ref{ahlfors-condition}.}  We apply
Lemma~\ref{polynomial-growth-implies-robustness} as follows: If $r <
R/2$, this lemma gives an upper bound on the $\lambda$-dimensional
Hausdorff measure of any ball $\opr{B}(x, r/2)$ for $r \leq R/2$ of
$2^{3 \, \lambda} \, k^{-1} \, K \, r^\lambda$.

Consider the case $E$ is open with $\nu(E) < \infty$. Let $V_i
\subseteq E$ be a disjoint sequence of closed balls of radius $r_i
\leq R/2$ such that either $\sum_i r_i^\lambda = +\infty$ or
$\mathcal{H}^\lambda(E \setminus \bigcup V_i) = 0$.  In the first
case,
$$\nu(E) \geq \sum_i \nu(V_i) \geq k \, \sum_i
r_i^\lambda = +\infty,$$
which contradicts $\nu(E) < \infty$. Thus for some constants $K_1,
K_2, K_3 $,
$$\mathcal{H}^\lambda(E) = \sum_i \mathcal{H}^\lambda(V_i) \leq \sum_i
K_1 \, r_i^\lambda \leq K_2 \sum_i \nu(V_i) \leq
K_3 \, \nu(E).$$
Open balls of radius $\leq R/2$ have finite $\nu$ measure (by
assumption) and finite $\mathcal{H}^\lambda$ measure by the preceding
lemmas. By the first part of the theorem, balls of radius $\leq R/2$
are relatively compact, hence separable.  By approximation properties
of finite Borel measures (Theorem 1.1 of~\cite{billingsley1}), it
follows that for open balls $V$ of radius $\leq R/2$
inequality~(\ref{counting-upper-bound-for-hausdorff-measure}) holds
for arbitrary Borel sets $E \subseteq V$.  Since each separable Borel
set $E$ is a countable union of Borel sets contained in balls of
radius $\leq R/2$, the general case follows.~$\blacksquare$

\begin{definition}\label{ahlfors-regularity-definition}
A complete extended metric space for which there are $0<k \leq K <
\infty$ and $R \in ]0, \infty]$ such that
\begin{equation}
\label{ahlfors-hausdorff-measure-condition} k \, r^\lambda \leq
\mathcal{H}^\lambda \bigl(\opr{B}(x,r)\bigr) \leq K \, r^\lambda
\end{equation}
for $0<r < R$ is called Ahlfors regular of dimension $\lambda$ up to
$R$.
\end{definition}
We emphasize that this definition allows $R = +\infty$.  Note that
Ahlfors regularity up to $+\infty$ is a property about the large scale
structure of a space unlike its dimensional behavior which is a local
property.

See~\cite{davidsemmes} for more on this circle of ideas.

\subsection{Discrete Characterization of Ahlfors Regularity}

We can easily obtain a converse to
Proposition~\ref{measure-characterization-polynomial-growth}.

\begin{prop}\label{non-standard-characterization-of-Ahlfors-regular}
A necessary condition a complete metric space $(X,d)$ be Ahlfors
regular of dimension $\lambda$ up to some $R \in ]0, \infty]$ is that
\begin{enumerate}
\item \label{item-1-ahlfors-characterization} There exist an internal
metric space $(\mathfrak{X}, d_{\mathfrak{X}})$ of polynomial growth
$\lambda$ on $[R_{-}, R_{+}]$ where $R_{-} \cong 0$ and $0 \ll R_{+}
\ll \infty$ if $R$ is finite or $R_{+}$ is an unlimited hyperreal in
case $R = \infty$
\item \label{item-2-ahlfors-characterization} $(X,d)$ is a subset of
$\lhull(\mathfrak{X},d_{\mathfrak{X}})$ with the property that for all
$r < R$ and all $x \in X$, $B(x,r) \subseteq X$.
\end{enumerate}
\end{prop}

\begin{lemma}
Suppose $\mu$ satisfies
inequality~(\ref{ahlfors-regular-measure-condition}) for all $r \in
]0, R[$.  Then there are constants $k', K'$ such that for all $\sigma
> 0$ and $Y \subseteq X$ which is maximal $\sigma$-separated,   the
following holds:
\begin{equation}\label{approximation-bound-standard}
k' \leq \frac{\opr{card}\bigl(\opr{B}(a, r)\cap Y \bigr)}{
{\bigl(\frac{r}{\sigma}\bigr)}^\lambda} \leq K'
%
%
\quad \mbox{if $ 2 \, \sigma \leq r < R/2$ and $a \in X$.}
\end{equation}
\end{lemma}
{\sc Proof.}  Let
$D=D(a, r, \sigma, Y)$ be the expression in
inequality~(\ref{approximation-bound-standard}) and let $k$ and $K$ be
as in inequality~(\ref{ahlfors-regular-measure-condition}).
Since $Y$ is $\sigma$-separated, the balls $\opr{B}(z,
\sigma/2)$ for $z \in Y$ are pairwise disjoint.
Thus
$$
\mu \bigl(\opr{B}(a, 2 \, r )\bigr) \geq \sum_{z \in \opr{B}(a, r)
\cap Y} \mu \bigl(\opr{B}\bigl(z, \frac{\sigma}{2}\bigr)\bigr)
\geq D \, \biggl(\frac{r}{\sigma}\biggr)^\lambda \, k \,
\biggl(\frac{\sigma}{2}\biggr)^\lambda 
=  2^{-\lambda} D \, k \, r^\lambda.
$$
Since $2 \, r < R$,
inequality~(\ref{ahlfors-regular-measure-condition}) and the previous
inequality imply $D \leq 4^\lambda \, K \, k^{-1}$.  In the other
direction,
$\opr{B}\bigl(a, \fraco{r}{2}\bigr) \subseteq \bigcup_{z \in \opr{B}(a,
r) \cap Y} \opr{B}(z, \sigma)$.
Thus,
$$\mu \bigl(\opr{B}\bigl(a, \frac{r}{2}\bigr)\bigr) \leq \sum_{z \in \opr{B}(a, r)
\cap Y} \mu\bigl( \opr{B}(z, \sigma)\bigr) \leq D \,
\biggl(\frac{r}{\sigma}\biggr)^\lambda \, K \, \sigma^\lambda.$$
Since $r/2 > 0$, inequality~(\ref{altrenate-polynomial-growth-equation})
again implies $D(a, r, \sigma, Y) \geq 2^{-\lambda} \, K^{-1} \,
k$.~$\blacksquare$

To complete the proof of
Proposition~\ref{non-standard-characterization-of-Ahlfors-regular},
for each real number $\sigma > 0$, let $Y_\sigma$ be a maximal
$\sigma$-separated subset of $X$.  Consider $\transfer{}$ of the
function $\sigma \rightarrow Y_\sigma$.  By transfer,
$[\transfer{Y}]_\sigma \subseteq \transfer X$
satisfies~(\ref{approximation-bound-standard}) and
$[\transfer{Y}]_\sigma$ is maximal $\sigma$ separated in $\transfer X$
for every positive hyperreal $\sigma$.  Taking $\sigma \cong 0$, it
follows $\phull{[\transfer{Y}]_\sigma} = \phull{\transfer{X}}$.  Now
let $\mathfrak{X} = Y_\sigma$. $\mathfrak{X}$
satisfies~(\ref{item-1-ahlfors-characterization}) with $R_{-} = 2 \,
\sigma$ and $R_{+}=R/2$.
To prove~(\ref{item-2-ahlfors-characterization}), by Ahlfors
regularity up to $R$ of $X$, for $s < R$ the closed balls
$\opr{\overline{B}}(x,s)$ in $X$ are compact and so
$\opr{\overline{B}}_{\phull{\transfer{X}}}(x,s) =
\opr{\overline{B}}(x,s) \subseteq X$.~$\blacksquare$

The characterization of Ahlfors regular spaces given by
Proposition~\ref{non-standard-characterization-of-Ahlfors-regular}
provides no information about the space
$(\mathfrak{X},d_{\mathfrak{X}})$.  The remainder of the paper is
devoted to obtaining a more informative result on the form of
$(\mathfrak{X},d_{\mathfrak{X}})$.

\section{Graph Spaces}

Graph means undirected graph, {\em without looping edges}.  We denote
the adjacency relation between vertices $x, y$ in $G$ by $x
\longleftrightarrow y$.  If $G$ is connected, the {\em graph distance}
between $x,y \in G$ is the length of the shortest path in $G$ from $x$
to $y$.  The graph distance is a metric with values in the nonnegative
integers; we denote the graph distance by $\opr{dist}$. For each $\Rescale
\in \mathbb{R}$, let $\opr{dist}_\Rescale(x,y) =
\frac{1}{\Rescale}\opr{dist}(x,y)$. We use $\opr{B}$ to denote open balls
relative to $\opr{dist}$ and $\opr{B}_\Rescale$ to denote open balls relative
to $\opr{dist}_\Rescale$.

We will be mainly considering the $\transfer$-versions of
graph-theoretic concepts.  Vertices $x,y$ of $G$ belong to the same
limited component of $(G, \opr{dist}_\Rescale)$ iff there is a path
from $x$ to $y$ of length $\preceq_O \Rescale$.  In particular, if the
graph diameter of $(G, \opr{dist}_\Rescale)$ is $\preceq_O \Rescale$,
then $(G, \opr{dist}_\Rescale)$ has exactly one limited component.

We contrast our notion of polynomial growth in metric spaces with the
standard one for graphs: If $G$ is a graph, $G$ has standard
polynomial growth $\lambda \in \mathbb{R}$ iff there are constants $0
< c \leq C < \infty$ in $\mathbb{R}$ such that
\begin{equation}\label{standard-polynomial-growth-inequality}
c \ \ell^\lambda \leq \opr{card} \opr{B}(x, \ell) \leq C
\ell^\lambda \quad \mbox{
for $\ell \geq 1$.}
\end{equation}
Note that standard polynomial growth for internal graphs is an
external property.

%
\begin{example} There is a vast literature on growth of finitely
generated discrete groups leading to Gromov's theorem on virtually
nilpotent groups (See~\cite{bass}, \cite{gromov}, \cite{gromov-1},
\cite{dries-wilkie}.) See~\cite{hall} for definitions and examples of
nilpotent groups.

We will need a basic result on nilpotent groups due to
Bass~\cite{bass}.  Recall first that for a group $G$, and $V \subseteq
G$ such that $V= V^{-1}$ and $e \not\in V$, the {\em Cayley graph} of
$G$ relative to $V$, denoted $\opr{Cayley}(G,V)$, is the undirected
loopless graph whose vertices are the elements of $G$ and whose edges
are the pairs $\{x,y\}$ such that $x^{-1} \, y \in V$.
%
%
If $G$ is a finitely generated nilpotent group, then for any symmetric
generating set $W$, there exists a $d(G) \in \mathbb{N}$ such that
$\opr{Cayley}(G,W)$ has standard growth $d(G)$.
%
%
Theorem 2 of~\cite{bass}, provides an explicit formula for $d(G)$.

\begin{prop}
Suppose $G$ is an internal nilpotent group of class $n \in
\mathbb{N}$.  Suppose furthermore $W \subseteq G$ is a symmetric set
with $\opr{card} W \in \mathbb{N}$ and such that $W$ generates $G$ as
an internal group. Then there is a standard integer $\lambda$ such
that $\opr{Cayley}(G,W)$ is of polynomial growth $\lambda$ on some
interval $[2, M]$ for $M \cong \infty$.
\end{prop}
{\sc Proof.}  The external group $G_\infty$ generated by $W$ is
nilpotent of class $\leq n$ and finitely generated. By Bass's result,
$\opr{Cayley}(G_\infty,W)$ is of standard polynomial growth $\lambda
\in \mathbb{N}$.  This means that for some $c,C \in ]0,\infty[$,
\begin{equation}\label{modified-standard-polynomial-growth-inequality}
c \, \ell^\lambda \leq \opr{card} \bigl(\opr{B}(e, \ell) \cap
G_\infty \bigr) \leq C \, \ell^\lambda
\end{equation}
for all $\ell \in \mathbb{N}$.  For $\ell \in \mathbb{N}$, 
$\opr{B}(e, \ell) \cap G_\infty = \opr{B}(e, \ell)$.
Thus~(\ref{modified-standard-polynomial-growth-inequality}) actually
implies an internal condition,
i.e. inequality~(\ref{standard-polynomial-growth-inequality}) with $x
= e$.  Thus by overspill there is an $R_{+} \cong \infty$ such
that~(\ref{standard-polynomial-growth-inequality}) holds with $x = e$
for all $\ell \leq R_{+}$.
However, $\opr{card} \opr{B}(e, \ell) = \opr{card} \opr{B}(x, \ell)$
for any $x \in G$, so $\opr{Cayley}(G,W)$ is of polynomial growth
$\lambda$ on some interval $[2, M]$ for $M \cong
\infty$.~$\blacksquare$
\end{example}

\section{Lifting Measures}

\begin{prop} 
\label{measlift-unbounded} Let $(X,d_X)$ be a $\sigma$-compact metric
space, $\mu$ a ball Borel measure on $X$ such that $\mu(K) < \infty$
for every compact $K\subseteq X$, $(\mathfrak{X},d_\mathfrak{X})$ a
hyperfinite space such that $\lhull{\mathfrak{X}} \supseteq X$.  Then
there is a hypermeasure $\nu$ on $\mathfrak{X}$ such that
$\opr{P}_\mathfrak{X}(\nu) \, | \, \lhull{\mathfrak{X}} = \mu$.
\end{prop}
{\sc Proof.}  This is a special case of the discussion
in~\cite{albeverio-etal}, \S5.2 on lifting measures.  The arguments
used there rely on more saturation than we need in this paper, but
given the separability assumptions on $(X,d_X)$ we can get by with
only countable saturation.  For completeness we give a direct proof
assuming only countable saturation.

Assume first $X$ is compact.  We may assume without loss of generality
that $\mu$ is a probability measure and $\mathfrak{X} = \transfer{X}$.
To show this last remark, note there is a hyperfinite $F \subseteq
\mathfrak{X}$ such that $\lhull{F} = X$ and a mapping $g: \transfer{X}
\rightarrow F$ which makes the following diagram commutative:
\begin{diagram} \transfer{X}  & \rTo^{g} & F \\
\dTo<{\varphi_{\transfer{X}}} & & \dTo>{\varphi_{F}} \\
X & \rTo{\opr{id}_X} & X
\end{diagram} 
Thus any lifting for $\mu$ in $\transfer{X}$ can be pushed over by $g$
to a lifting for $\mu$ in $F$.
There is a countable set $\mathcal{K}$ of compact sets in $X$ such
that every open set in $X$ is the union of a nondecreasing sequence of
$\mathcal{K}$.  By the saturation property, there is a hyperfinite
boolean algebra $\mathcal{F} \subseteq \transfer \mathfrak{B}$, where
$\mathfrak{B}$ is the Borel algebra of the compact metric space $X$,
such that for each $K \in \mathcal{K}$, $\transfer K \in \mathcal{F}$.
Now there is an internal operator $\kappa$ on the limited nonnegative
hyperfinitely additive set functions on $\mathcal{F}$ such that for
any limited hyperfinitely additive set function $\nu$, $\kappa(\nu)$
extends $\nu$ and is a limited hypermeasure on $\transfer X$.
$\transfer \mu$ is a limited nonnegative hyperfinitely additive set
function on the algebra $\transfer \mathfrak{B}$.  Let $\nu=\kappa
\left( \transfer \mu | \mathcal{F}\right)$ and $\rho = \opr{P}(\nu)$.
$\transfer \mu$ has total mass one and thus $\nu$ has total mass
one. $\rho$ is a Borel measure on $X$ of total mass $\leq 1$.  For any
$K \in \mathcal{K}$, $\transfer K \subseteq \varphi^{-1}(K)$ so that
$$\rho(K)={\opr{L}(\nu)}(\varphi^{-1}(K)) \geq
{\opr{L}(\nu)}(\transfer K) = \opr{st} (\transfer \mu(\transfer K)) =
\mu(K).$$ 
By monotonicity of measures, and the fact every open set is the union
of a nondecreasing sequence of $\mathcal{K}$, it follows that $\rho(U)
\geq \mu(U)$ for any open set $U$.  By regularity
(see~\cite{billingsley1}, Theorem 1.4) , $\rho(E) \geq \mu(E)$ for any
Borel set $E$.  Since $\mu$ is a probability measure and $\rho(X) \leq
1$, $\rho$ is also a probability measure.  It follows $\rho = \mu$.

In the general case, there is a nondecreasing sequence of compact sets
$\{K_i\}_{i \in \mathbb{N}}$ such that $X= \bigcup_i K_i$.  For each
$i \in \mathbb{N}$, let $E_i \subseteq \mathfrak{X}$ be hyperfinite
and such that $\lhull{E_i} \supseteq K_i$.  Existence of the sets
$E_i$ follows immediately from countable saturation.  Clearly we may
assume $E_i \subseteq E_{i+1}$.  By assumption, $\mu | K_i$ is a
finite Borel measure.  By the special case considered in the previous
paragraph, there is a sequence of limited hypermeasures $\{\nu_i\}_{i
\in \mathbb{N}}$ such that $\nu_i$ is supported in $E_i$ and
$\opr{P}(\nu_i) = \mu | K_i$.  We may assume $\nu_i \leq
\nu_{i+1}$. Proof: Consider $\mu | K_i$ as a measure on $K_{i+1}$ by
assigning total measure $0$ to $K_{i+1} \setminus K_i$.  With this
convention, $\mu | K_i \leq \mu | K_{i+1}$, so that $\mu | K_{i+1} -
\mu | K_{i}$ is also a finite countably additive measure.  Now let
$\nu'_{i+1}$ be such that $\opr{P}(\nu_i') = \mu | K_{i+1} - \mu |
K_{i}$ and $\nu'_{i+1}$ is supported in $E_{i+1}$.  Then $\nu_{i+1} =
\nu'_{i+1} + \nu_i$ is the desired measure.

$\varphi^{-1}(K_i)$ is Borel, so by approximation of Loeb measurable
sets of finite measure by internal subsets, there is an internal set
$F_i$ such that $\varphi^{-1}(K_i) \subseteq F_i$ and
$$\nu_{i+1}(F_i) \leq \opr{L}(\nu_{i+1})\bigl(\varphi^{-1}(K_i)\bigr)
+ 2^{-(i+1)}.$$
Now 
$$\opr{L}(\nu_{i+1})\bigl(\varphi^{-1}(K_i)\bigr) =
\opr{P}(\nu_{i+1})(K_i) = \mu(K_i) = \opr{P}(\nu_{i})(K_i) =
\opr{L}(\nu_{i})\bigl(\varphi^{-1}(K_i)\bigr), $$
so
$$
\nu_{i+1}(F_i) \leq \opr{L}(\nu_i)\bigl(\varphi^{-1}(K_i)\bigr) +
2^{-(i+1)} \leq \opr{L}(\nu_i)\bigl(F_i) + 2^{-(i+1)}
\leq  \nu_i(F_i) + 2^{-i}.
$$
Since $\nu_i \leq \nu_{i+1}$ it follows that for all internal $A
\subseteq F_i$,
$$\nu_i(A) \leq \nu_{i+1}(A) \leq \nu_i(A) + 2^{-i}.$$

Now apply saturation and overspill, to conclude that the sequences
$\nu_k$, $E_k$ and $F_k$ have extensions to hyperfinite
sequences defined for $1 \leq k \leq N$ such that
\begin{enumerate}
\item $E_i \subseteq E_{i+1}$ for $i \leq N-1$,
\item $E_i$ is hyperfinite and $E_i \subseteq F_i$ for $i \leq N$,
\item \label{nesting-of-measures} $\nu_i(A) \leq \nu_{i+1}(A) \leq
\nu_{i}(A) + 2^{-i}$ for $i \leq N-1$ and all internal $A
\subseteq F_i$.
\end{enumerate}
It follows from item~\ref{nesting-of-measures}, that
for $j \geq i$,
\begin{equation}
\label{nesting-of-measures-1}\nu_i(A) \leq \nu_j(A) \leq \nu_{i}(A) + 2^{-(i-1)}
\end{equation}
for all internal $A \subseteq F_i$.  It follows immediately
from~(\ref{nesting-of-measures-1}) and the monotone class theorem that
for any inner Borel set $A \subseteq F_i$,
$$\opr{L}(\nu_i)(A) \leq \opr{L}(\nu_j)(A) \leq \opr{L}(\nu_i)(A) +
2^{-(i-1)}.$$

Let $\nu = \nu_N$.  We need to show that for any Borel set $A
\subseteq \lhull{\mathfrak{X}}$, $\opr{P}(\nu)(A) = \mu(A)$.  If $i
\in \mathbb{N}$, then for any ball Borel set $A \subseteq K_i$,
$\varphi^{-1}(A) \subseteq \varphi^{-1}(K_i) \subseteq F_i$, so
\begin{eqnarray*}
\mu(A) = \opr{P}(\nu_i)(A) & = &
\opr{L}(\nu_i)\bigl(\varphi^{-1}(A)\bigr) \\
& \leq & \opr{L}(\nu)\bigl(\varphi^{-1}(A)\bigr) \\
& = & \opr{P}(\nu)(A) \\
& \leq &  \opr{L}(\nu_i)\bigl(\varphi^{-1}(A)\bigr) + 2^{-(i-1)} \\
& = & \opr{P}(\nu_i)(A) + 2^{-(i-1)} = \mu(A) + 2^{-(i-1)}
\end{eqnarray*}
Since $i$ is arbitrary, $\mu(A) = \opr{P}(\nu)(A)$ for any Borel set
$A$ which is a subset of some $K_i$. However, $\bigcup_{i \in
\mathbb{N}} K_i = X$, so by countable additivity of $\mu$ and of
$\opr{P}(\nu)$, the result follows for all Borel sets $A \subseteq
X$.~$\blacksquare$


A hypermeasure on the internal subsets of $\mathfrak{X}$ is {\em
uniform} if all the singleton sets have the same mass.  Define the
{\em density} of $\nu$ to be the mass of each singleton.

\begin{prop} \label{even-covering}
Let $\nu$ be a hypermeasure on a hyperfinite set $\mathfrak{X}$. Then there is a
uniform hypermeasure $\mu$ on a hyperfinite set $\mathcal{Y}$ and an
internal mapping $p:\mathcal{Y} \rightarrow \mathfrak{X}$ such that $\mu(p^{-1}(A)) \cong
\nu(A)$ for every internal $A \subseteq \mathfrak{X}$ and $\mu(\mathcal{Y}) \leq
\nu(\mathfrak{X})$.
\end{prop}
{\sc Proof.}
Assume $\mathfrak{X}=\{1, \ldots, N\}$ and let $\nu_1, \ldots, \nu_N$ be the
masses of the atoms of $\mathfrak{X}$. Let $M$ be such that $M/N \cong \infty$.
For each $i \in \{1, \ldots, N\}$ there is a unique hyperinteger $0
\leq k_i$ such that
$k_i/M \leq \nu_i < (k_i+1)/M$.
Let $\{B_i\}_{1 \leq i \leq N}$ be disjoint sets such that
$\opr{card}(B_i)=k_i$, and define $p:\mathcal{Y} \rightarrow
\mathfrak{X}$ so that $p$ is identically $i$ on $B_i$, for $i \geq
1$. Let $\mu$ be the hypermeasure on $\mathcal{Y}$ which assigns the
mass $1/M$ to each singleton of $\mathcal{Y}$. For every internal $A
\subseteq \mathfrak{X}$,
$$ \mu(p^{-1}(A)) = \sum_{i \in A} \mu(B_i) = \sum_{i \in A} k_i/M =
\sum_{i \in A}(\nu_i+o_i) = \nu(A)+\sum_{i \in A} o_i $$
Now $|\sum_{i \in A} o_i | \leq \sum_{i \in \mathfrak{X}} | o_i |\leq
N/M \cong 0$.  Finally $\mu(\mathcal{Y}) = \sum_i k_i/M \leq \sum_i
\nu_i = \nu(\mathfrak{X})$.~$\blacksquare$

It is clear that there is a lot of leeway in the choice of density
$1/M$.  In particular, $M$ can be chosen to be a hyperinteger.
In case $\mu$ is a hypermeasure on $\mathfrak{X}$ and
$(\mathfrak{X},d)$ is an internal metric space, we can take
$\mathcal{Y}$ to be an internal metric space.  In case
$(\mathfrak{X},d)$ is a graph metric space we can $\mathcal{Y}$ to be
a graph metric space also.

\begin{definition}\label{definition-of-spiked-graph}
A {\em spiked graph over $G$} is a family of connected loopless graphs
$\{H_v\}_{v \in \opr{nodes}(G)}$ which are pairwise disjoint and such
that $v \in H_v$ for each node $v$ of $G$.

The {\em disjoint union} of a spiked graph is the loopless graph $G'$ obtained
as follows: The set of nodes of $G'$ is the disjoint union
$\bigsqcup_v \opr{nodes}(H_v)$. An edge of $G'$ is either an edge in
$G$ or an edge in one of the graphs $H_v$. We will use the expression 
$\bigsqcup_v H_v$ to denote the disjoint union.

The {\em base projection} of the disjoint union is the mapping which
is identically $v$ on $H_v$.
\end{definition}

We will use informal but suggestive terminology to describe spiked
graphs: The graphs $H_v$ are the spikes, etc.. Note that the disjoint
union of a spiked graph over $G$ contains $G$ as a full subgraph.  A
{\em spike extension} of a graph $G$ is a loopless graph $G'$ which is the
disjoint union of a spiked graph over $G$.

\begin{prop} \label{even-covering-metric-case}
Under the assumptions of Proposition~\ref{even-covering}, if $(\mathfrak{X},d_\mathfrak{X})$
is an internal metric space, we may assume $(\mathcal{Y},d_\mathcal{Y})$ is an internal
metric space and $p: \mathcal{Y} \rightarrow \mathfrak{X}$ satisfies
$$d_\mathfrak{X}(p(y),p(y')) \leq d_\mathcal{Y}(y,y') \leq d_\mathfrak{X}(p(y),p(y')) +
\epsilon(y,y')$$
where $\epsilon(y,y')\cong 0$.  

Suppose $d_\mathfrak{X}$ is $\delta \times \opr{dist}$ where $\delta \cong 0$ and
$\opr{dist}$ the metric associated to a loopless graph on $\mathfrak{X}$.  Then $d_\mathcal{Y}$ can
also be taken as $\delta \, \opr{dist}'$ for some metric associated to
a loopless graph on $\mathcal{Y}$ which is a spike extension of $\mathfrak{X}$ and $p$ is the base
projection.
\end{prop}
{\sc Proof.} Provide each set $p^{-1}(x)$ with a metric $d'_x$ which
assumes only infinitesimal values and with a distinguished point
$y_x$.  The metric $d'_x$ can be chosen in a completely arbitrary way,
but one possible choice is as follows: Choose an infinitesimal
$\theta$ and $d'_x(a,b) = \theta$, whenever $a \neq b$.  Define
$$d_\mathcal{Y}(y,y') = d_\mathfrak{X}(p(y), p(y'))+d'_{p(y)}(y, y_{p(y)})+d'_{p(y')}(y_{p(y')}, y')$$
It is clear $d_\mathcal{Y}$ has the desired properties.

In the case $d_\mathfrak{X} = \delta \times \mbox{graph metric on
$\mathfrak{X}$}$, consider $\mathcal{Y}$ as a loopless graph obtained
by attaching spikes to $\mathfrak{X}$, where now each spike is a
complete graph. Then $p : \mathcal{Y} \rightarrow \mathfrak{X}$ is the
mapping which maps each member of $\mathcal{Y}$ to the attachment
point of the spike in $\mathfrak{X}$.  Since each spike is a complete
graph, $\opr{dist}(y,y') \leq 1$ for all $y,y'$ on a spike, the metric
$d_\mathcal{Y}$ on each spike is infinitesimal.~$\blacksquare$

The space $\mathcal{Y}$ of the preceding result can be regarded as an
infinitesimal thickening of $\mathfrak{X}$.  Note that also the
following diagram commutes 
\begin{diagram} \mathcal{Y} & \rTo^{p} &
\mathfrak{X}\\
\dTo<{\varphi_{\mathcal{Y}}} &  &
\dTo>{\varphi_{\mathfrak{X}}} \\
\lhull{\mathcal{Y}} & \rTo & \lhull{\mathfrak{X}}
\end{diagram}
where the bottom arrow is an isometric isomorphism.

\subsection{Graph Regularity}

Nothing has been said about a {\em regular} graph structure on
$\mathcal{Y}$ (e.g., all nodes have the same degree), but this can
also be arranged.

\begin{prop} \label{even-covering-metric-case-regular}
Under the assumptions of Proposition~\ref{even-covering}, suppose
$d_\mathfrak{X}$ is $\delta \times \opr{dist}$ where $\delta \cong 0$
and $\opr{dist}$ is the metric associated to a loopless graph on
$\mathfrak{X}$.  Then $d_\mathcal{Y}$ can also be taken as $\delta \,
\opr{dist}'$ for some metric associated to a regular loopless graph on
$\mathcal{Y}$ which is a spike extension of $\mathfrak{X}$.
\end{prop}
{\sc Proof.} By Proposition~\ref{even-covering-metric-case}, we may
assume $\nu$ is a uniform hypermeasure on $\mathfrak{X}$ with density
$M$.  Let $k \in \transfer{\mathbb{N}}$ be odd and such that $k \geq
3+\max_{x \in \mathfrak{X}} \opr{deg}_{\mathfrak{X}}(x)$.  Since
$\mathfrak{X}$ is hyperfinite, such a $k$ exists.  There exists a
spiked graph over $\mathfrak{X}$, $\{F_x\}_{x \in \mathfrak{X}}$ such
that
\begin{enumerate}
\item $\opr{card} F_x = k+3$,
\item For each $x \in \mathfrak{X}$, $x$ has degree $k + 1 -
\opr{deg}_{\mathfrak{X}}(x)$ in $F_x$,
\item For each $x \in \mathfrak{X}$ and all $y \in F_x \setminus
\{x\}$, $y$ has degree $k + 1$ in $F_x$.
\end{enumerate}
Then the loopless graph $\mathcal{Y}=\bigsqcup_{x \in \mathfrak{X}}
F_x$ (see Definition~\ref{definition-of-spiked-graph}) is regular of
degree $k + 1$.  If $\opr{dist}'$ is the graph distance on
$\mathcal{Y}$, then each $x \in \mathfrak{X}$ is at $\opr{dist}'$
distance at most $2$ from every $y \in F_x$. To see this, note that by
the degree assumption on $x$, $x$ is connected to at least one $y_x
\in F_x \setminus \{x\}$. By the degree assumption on $y$, every $y
\in F_x \setminus \{x\}$ is connected to all other $y' \in F_x$ but
one. If $y$ is connected to $x$, we are done. Otherwise, $y$ is
connected to $y_x$.
It follows that each $x \in \mathfrak{X}$ is at infinitesimal
$d_{\mathcal{Y}}$ distance from every $y \in F_x$. Let $\nu'$ be the
uniform hypermeasure on $\mathcal{Y}$ with density $M/(k+3)$.

We need to prove that for each $x$ there is a loopless graph $F_x$
with the required property.  Start off with a disjoint family of sets
$\{F_x\}_{x \in \mathfrak{X}}$ such that $x \in F_x$ and $\opr{card}
F_x = k+3$.  We will add edges to $F_x$ so that the above requirements
are met.  Partition $F_x$ into the sets $\{x\}$, $A$ and $B$ as shown
in the figure, where $s = k+1 - \opr{deg}_{\mathfrak{X}}(x) \geq 4$:
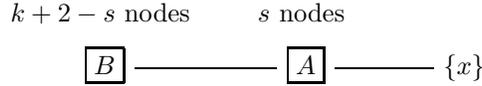
\begin{diagram}
k+2-s \text{ nodes }  &   &  s \text{ nodes }  &  & \\
\boxed{B} &  \rLine & \boxed{A} &  \rLine & \{x\} 
\end{diagram}
The line on the right represents a set of $s$ edges joining $x$ to
each of the $s$ members of $A$. As the following lemma shows, we can
add edges between members of $A$ so that each element of $A$ is joined
to exactly $s-2$ other elements of $A$:

\begin{lemma}
Suppose $\opr{card} A \geq 4$ is even and $r \leq \opr{card} A- 2$ is
even.  Then there is a connected loopless graph on $A$ which is regular of
degree $r$.
\end{lemma}
{\sc Proof.} Let $s = \opr{card} A$ and consider the additive group
$\mathbb{Z}/(s)$.  $\mathbb{Z}/(s)$ is a cyclic group of order
$s$. Let $V$ be the subset of $\mathbb{Z}/(s)$ consisting of the
equivalence classes of $\pm 1, \ldots, \pm r/2$.  $V$ has cardinality
$r$. The Cayley graph $\opr{Cayley}(\mathbb{Z}/(s), V)$ is regular of
degree $r$.  Since $0 \not\in V$, $\opr{Cayley}(\mathbb{Z}/(s), V)$ is
also loopless.~$\blacksquare$

Each element of $A$ so far has attached $s-1$ edges : $s-2$ connected
to other members of $A$ and $1$ edge to $x$.  To get up to $k+1$ edges
per node, we must leave $k+1-(s-1) = k+2-s$ unattached out edges for
each of the $s$ members of $A$.  This is a total of $s ( k+2-s)$
dangling edges (to the left) from $A$.  Similarly each element of $B$
can be joined to $k+1-s$ nodes of $B$ --- just consider the complete
graph on $B$--- leaving $k+1- (k+1-s) = s$ unattached out edges per
element of $B$. This is a total of $s (k+2 -s)$ dangling edges (to the
right) from $B$. Now attach each right edge from $B$ to a unique left
edge from $A$.~$\blacksquare$

\section{Spaces from Graphs}
\label{homogeneous-space-presentation}

\begin{example}
Consider $\mathbb{R}^n$ with the metric $d(\vec{x}, \vec{y}) = \sum_k
| x_i - y_i|$.  This metric is sometimes called the Taxicab
distance~\cite{taxicab}.  Closed $n$-cubes $\prod_{i=1}^n [a_i, b_i]$
in $\mathbb{R}^n$ have the length space property with respect to the
Taxicab metric.  \end{example}

\begin{example}
More generally, in the case $n \geq 2$ it is geometrically clear that
the result of removing a finite number of open $n$-cubes from a
closed $n$-cube is also a length space.  This example easily extends
to the removal of a countable number of disjoint open $n$-cubes.
\end{example}

\begin{example}
Suppose $G$ is a hyperfinite loopless graph with graph distance $\opr{dist}$.
For any $M \cong \infty$, $\lhull{(G, \opr{dist}_M)}$ is a length
space.  Proof: Suppose $\hat{x}, \hat{y} \in \lhull{G}$ are such that
$\opr{st} \bigl(\opr{dist}_M(x,y)\bigr) =
\lhull{\opr{dist}_M}(\hat{x}, \hat{y}) = a$.  Thus there is a path $x
= z_0 \longleftrightarrow z_1 \cdots \longleftrightarrow z_n = y$ with
$n = \opr{dist}(x,y)$ and
$a = \opr{st} (n / M)$.  The map $T = \{0, 1/M, \ldots, n/M\}
\rightarrow (G, \opr{dist}_M)$ given by $j/M \mapsto z_j$ is clearly
isometric and thus factors through an isometry $[0, a] \rightarrow
\lhull{(G, \opr{dist}_M)}$.

\end{example}

In this section we will characterize Ahlfors regular length spaces
$(X,d_X)$ as bounded components of hulls of hyperfinite spaces
$(\mathfrak{X},d)$.  The main technical point of the proof is to show
that a certain subset $\mathcal{Y}$ of $\mathfrak{X}$ has polynomial
growth.  Specifically, we need estimates of the kind
$$k \, r^\lambda \leq \opr{card}\bigl(\opr{B}(x,r) \cap \mathcal{Y} \bigr) \leq
K \, r^\lambda$$
for $x \in \mathcal{Y}$. Note that this lower bound is stronger than
the lower bound
$k \, r^\lambda \leq \opr{card}\opr{B}(x,r)$.
%

\begin{lemma}
Suppose $(X,d)$ is a length space and $\rho > 0$. Let $a, x \in X$ be
such that $d(a,x) \leq \rho$.  Then for all $r$ such that $0 \leq r
\leq \rho$, there is a $c \in X$ such that
\begin{equation} \label{ball-nesting-condition-for-lenght-spaces}
\opr{B}\bigl(c, \frac{r}{2}\bigr) \subseteq \opr{B}(a, \rho) \cap \opr{B}(x, r).
\end{equation}
\end{lemma}
{\sc Proof.}  We consider two cases:

Suppose $d(a,x) < r/2$. In this case take $c = x$. Obviously,
$\opr{B}(x, \fraco{r}{2}) \subseteq \opr{B}(x, r)$. If $y \in
\opr{B}(x, r/2)$, then 
$d(a,y) \leq d(a,x)+d(x,y) < r/2 +r/2 = r \leq \rho$
so that $\opr{B}(x, r/2) \subseteq \opr{B}(a, \rho)$.

Suppose $d(a,x) \geq r/2$.  In this case we need the length space
property of $(X,d)$.  Let $\rho'= d(a,x)$.  There is an isometric map
$f: [0, \rho'] \rightarrow X$ with $f(0) = a$ and $f(\rho')= x$. Since
$0 \leq \rho' - r/2 \leq \rho'$, $c = f(\rho' - r/2)$ is
well-defined. Now $d(a,c) = \rho' -r/2$ and $d(c,x) = |\rho' - r/2 -
\rho'| = r/2$.  Thus by the triangle inequality,
$\opr{B}(c, \fraco{r}{2}) \subseteq \opr{B}(x,r)$
and
$\opr{B}(c, \fraco{r}{2})\subseteq \opr{B}(a, \rho') \subseteq
\opr{B}(a, \rho)$.~$\blacksquare$

\begin{lemma} \label{modified-lower-ahlfors-bound}
Suppose $(X,d)$ is an Ahlfors regular length space of dimension
$\lambda$ up to $R$.  Then there exists $k' > 0$ such for every $a,   \, x
\in X$ and $r, \, \rho > 0$ satisfying $r < R$ and $r \leq \rho$ and $
d(x,a) \leq \rho$,
$$ k' \, r^\lambda \leq \mathcal{H}^\lambda \bigl(\opr{B}(x,r) \cap
\opr{B}(a,\rho)\bigr)$$
\end{lemma}
{\sc Proof.} Let $k$ be as in
Formula~(\ref{ahlfors-hausdorff-measure-condition}).  By the preceding
lemma there is a $c$ satisfying the
inclusion~(\ref{ball-nesting-condition-for-lenght-spaces}).  Thus,
$$k \, \biggl(\frac{r}{2}\biggr)^\lambda \leq \mathcal{H}^\lambda
\bigl(\opr{B}(c,\frac{r}{2})\bigr) \leq \mathcal{H}^\lambda
\bigl(\opr{B}(x,r) \cap \opr{B}(a,\rho)\bigr)$$
so $k' = 2^{-\lambda} \, k$ will do.~$\blacksquare$

\begin{prop}\label{ahlfors-regular-characterization}
Suppose $(X,d_X)$ is an Ahlfors regular length space of dimension
$\lambda$ up to $R \in ]0, \infty]$, $(\mathfrak{X},d_\mathfrak{X})$
is a hyperfinite space such that
\begin{enumerate}
\item $(X,d_X)$ is a bounded component of
$\lhull{(\mathfrak{X},d_\mathfrak{X})}$,
\item For some uniform hypermeasure $\nu$ on $\mathfrak{X}$,
$\opr{P}_\mathfrak{X}(\nu) |_X = \mathcal{H}^\lambda$.
\end{enumerate}
Then there is an internal $\mathcal{Y} \subseteq \mathfrak{X}$ such
that $\lhull{\mathcal{Y}} \supseteq X$ and
$(\mathcal{Y},d_\mathfrak{X}| \mathcal{Y})$ is of polynomial growth
$\lambda$ on some interval $[R_{-} , \, R]$ with $R_{-} \cong 0$ in
case $R$ is finite or $[R_{-} , \, R_\infty ]$ for some unlimited
$R_\infty$ and $R_{-} \cong 0$ in case $R = \infty$.
\end{prop}
{\sc Proof.} Choose $a \in \mathfrak{X}$ so that $\hat{a} \in X$.  In
particular, $X$ is the bounded component of  $\hat{a}$.
%
%
%
For any standard $r$,
standard $\rho$ 
and $x \in \mathfrak{X}$ the
Formula~(\ref{ball-relation-between-beteen-space-and-hull}) implies
%
%
\begin{eqnarray*}
\opr{st} \, \nu\bigl(\opr{B}(x,\fraco{r}{2})\bigr) & = &
\opr{L}(\nu)\bigl(\opr{B}(x,\fraco{r}{2})\bigr) \\
& \leq & \opr{L}(\nu)\bigr(\varphi^{-1} \,
\opr{B}(\hat{x},r)\bigr) \\
& = & 
\opr{P}(\nu)\bigl(\opr{B}(\hat{x},r)\bigr) \\
& = & \mathcal{H}^\lambda\bigl(\opr{B}(\hat{x},r)\bigr)
\end{eqnarray*}
and similarly,
\begin{eqnarray*}
\opr{st} \nu\bigl(\opr{B}(x, r) \cap \opr{B}(a, \rho)\bigr) 
& = & \opr{L}(\nu)\bigl(\opr{B}(x, r) \cap \opr{B}(a, \rho)\bigr) \\
& \geq &  \opr{L}(\nu)\bigl(\varphi^{-1} \, \opr{B}(\hat{x}, r) \cap 
\varphi^{-1} \, \opr{B}(\hat{a}, \rho)\bigr) \\
& = & \opr{P}(\nu)\bigl(\opr{B}(\hat{x}, r) \cap \opr{B}(\hat{a},
\rho)\bigr) \\
& = & \mathcal{H}^\lambda\bigl(\opr{B}(\hat{x}, r) \cap \opr{B}(\hat{a},
\rho)\bigr).
\end{eqnarray*}
Now we use the bounds given by Ahlfors regularity with the modified
lower bound given by Lemma~\ref{modified-lower-ahlfors-bound}: If $x
\in \mathfrak{X}$ and $r, \rho$ are standard such that $r < R$, $r \leq \rho$ and
$d(a,x) \leq \rho$,
$$k \, r^\lambda \leq \mathcal{H}^\lambda\bigl(\opr{B}(\hat{x}, r)
\cap \opr{B}(\hat{a}, \rho)\bigr) $$
We translate these inequalities as follows.  Let $1/M$ be the density
of $\nu$.  If $x \in \mathfrak{X}$ and $r, \rho$ are standard such
that $r < R$ and $d(a,x) \leq \rho$,
$$\opr{st}\biggl(\frac{\opr{card} \opr{B}(x,r/2) }{M}\biggr)
\leq K \, r^\lambda$$
and if in addition $r \leq \rho$,
$$\opr{st}\biggl(\frac{\opr{card}\bigl(\opr{B}(x, r) \cap \opr{B}(a,
\rho)\bigr)}{M}\biggr) \geq k \, r^\lambda.$$
Thus, for standard $K' > K$ and $k' < k$ and all standard $\rho$ and
standard $r$ such that $0 < r < R$ and $x \in \mathfrak{X}$ such that
$d(a,x) \leq \rho$
$$\opr{card}\opr{B}(x,r/2) \leq M \, K' \, r^\lambda 
%
$$
and if in addition $r \leq \rho$,
$$\opr{card}\bigl(\opr{B}(x, r) \cap \opr{B}(a, \rho) \bigr) \geq M \,
k' \, r^\lambda 
%
%
$$

Now the same inequalities holds for arbitrary $r \in \transfer
\mathbb{R}$ for which $0 \ll r \ll R$ as is shown in the following
lemma:

\begin{lemma} \label{sufficient-sets-for-polynomial-growth}
Suppose $\psi$ is an internal nondecreasing function on $\transfer [0,
\infty[$ and $k, K$ are positive standard constants such that:
\begin{equation} \label{inequality-in-set-gamma-condition}
k \, s^\lambda \leq \psi(s) \leq K \, s^\lambda
\end{equation}
for all standard $s$ in some interval $]\alpha, \beta[$ with $\alpha,
\beta$ positive. Then
inequality~(\ref{inequality-in-set-gamma-condition}) holds for all $s
\in \transfer \mathbb{R}$ such that $\alpha \ll s \ll \beta$ with a
possibly different value for $K$ (not exceeding $(3/2)^\lambda$ the
original value of $K$) and for $k$ (but not less than $(2/3)^\lambda$
the original value of $k$).
\end{lemma}
{\sc Proof.} If $\alpha \ll s \ll \beta$, there is a standard $r \ll
\beta$ such that $s \leq r \leq 3/2 \, s$; for instance, take $r = 1/2
\, \bigl(\opr{st}(s) + \min (\opr{st}(\beta)+3 /2 \opr{st}(s))\bigr)$
in case $\beta$ is limited and $r = 5 / 4 \opr{st}(s)$ otherwise.
Thus
$$\psi(s) \leq \psi(r) \leq K \, r^\lambda \leq (3/2)^\lambda \, K \,
s^\lambda,$$
and similarly for the lower bound.~$\blacksquare$

\paragraph{Completion of Proof of Proposition
of~\ref{ahlfors-regular-characterization}.}  By the lemma, for
standard $\rho$, for $r \in \transfer \mathbb{R}$ such that $0 \ll r
\ll R$ and $x \in \mathfrak{X}$ such that $d(a,x) \leq \rho$
\begin{equation}\label{cardinality-growth-upper-bound}
\frac{\opr{card} \opr{B}(x, r)}{r^\lambda} \leq M \, K'
\end{equation}
and if $r \ll \rho$,
\begin{equation}\label{cardinality-growth-lower-bound}
\frac{\opr{card}\bigl(\opr{B}(x, r) \cap
\opr{B}(a, \rho) \bigr)}{r^\lambda}
\geq M \, k'
\end{equation}

For every $\rho \in \mathbb{N}$ let $r_\rho = 1/\rho$.  Then
\begin{enumerate}
\item Inequality~(\ref{cardinality-growth-upper-bound}) holds for $r$
such that $r_\rho \leq r \leq R - r_\rho$ and $d(a,x) \leq \rho$.  In
the case $R = +\infty$,
inequality~(\ref{cardinality-growth-upper-bound}) holds for limited $r$
such that $r_\rho \leq r$.
\item \label{insure-validity-of-inequality}
Inequality~(\ref{cardinality-growth-lower-bound}) holds for $r$ which
in addition satisfy $r \leq \rho- r_\rho$.
\end{enumerate}
In particular, by saturation and overspill there is a $\rho \cong
\infty$ and an $r_\infty \leq 1/\rho \cong 0$ such that these same
inequalities hold for limited $r$ such that $r_\infty \leq r \leq R -
r_\infty$.  Note that the additional restriction given by
item~(\ref{insure-validity-of-inequality}) to insure the validity of
inequality~(\ref{cardinality-growth-lower-bound}) disappears since $r$
is limited.  

Note that if $R=+\infty$, use overspill to conclude that these same
inequalities hold for $r$ such that $r_\infty \leq r \leq R_{+}$ for
some unlimited $R_{+}$.

To complete the proof, let $\mathcal{Y} = \opr{B}(a, \rho)$.~$\blacksquare$

\paragraph{Remarks.}  Note that the only property that we have used of
Hausdorff measure $\mathcal{H}^\lambda$ on $(X,d)$ is that it
satisfies an inequality of the
form~(\ref{ahlfors-regular-measure-condition}).  Thus the same result
is true of any Borel measure $\mu$ on $X$ which is boundedly
equivalent to $\lambda$-dimensional Hausdorff measure, that is which
is absolutely continuous with respect to Hausdorff measure and for
which the Radon-Nikodym derivative satisfies
$$c \leq \frac{d \mu}{d \mathcal{H}^\lambda}(x) \leq C \quad
\mbox{for almost all $x \in X$, }$$
where $0 < c \leq C < \infty$ and ``almost all'' is relative to
$\mathcal{H}^\lambda$.

\section{Representation of Ahlfors Regular Spaces}
Suppose $\delta$ is a hyperreal.  An internal metric space is
$\delta$-connected iff for every $x,y \in \mathfrak{X}$ with $\delta \leq d(x,
y)$, there is a sequence $x=x_0, \ldots , x_n=y$ such that
$d(x_i, x_{i+1}) \leq \delta$ for all $i \leq n-1$ and $n \delta \sim_O
d(x,y)$.
Gromov in~\cite{geometric-group-theory} introduces a related notion
called {\em long range connectedness}.  Note that our definition is
external.

\begin{example}
Suppose $(\mathfrak{X},d)$ is an internal length space.  Then for any $\delta >
0$, $(\mathfrak{X},d)$ is $\delta$-connected.  Proof: Let $\delta \leq a =
d(x,y)$.  By assumption, there is an internal isometric map $f: [0, a]
\rightarrow \mathfrak{X}$ such that $f(0)=x, f(a) = y$. Let $n$ be the largest
hyperinteger such that $(n - 1) \, \delta < a$ and define $t_k = k \,
\delta$ for $k \leq n-1$, $t_{n} = a$.  As succesive $t_k$'s differ by
less than $\delta$ and $f$ is an isometry, $ d(f(t_{k+1}), f(t_{k}))
\leq \delta$.  Moreover, $(n - 1) \, \delta < a \leq n \, \delta$ so
$1 \leq n \, \delta / a < 1 + \delta / a \leq 2$ and therefore $n \,
\delta \sim_O a$.  \end{example}

We will using the following internal function between subsets of
$\mathfrak{X}$: $d(A, B) = \sup_{x \in A} \inf_{y \in B } d(x,y)$.
This function, unlike the Hausdorff distance function on subsets is
not symmetric.
  
\begin{prop} \label{delta-subsets-prop}
Suppose $\mathcal{Y} \subseteq \mathfrak{X}$, and $d(\mathfrak{X},
\mathcal{Y}) \leq \delta$. If $\mathfrak{X}$ is $\delta$-connected,
then $\mathcal{Y}$ is $4 \delta$-connected.
\end{prop}
{\sc Proof.} Suppose $x,y \in \mathcal{Y}$ and $\delta < 4 \, \delta \leq d(x,
y)$.  By assumption there is a sequence $x=x_0, \ldots , x_n=y$ in $\mathfrak{X}$
such that $d(x_i, x_{i+1}) \leq \delta$ for all $0 \leq i \leq n-1$
and $n \delta \sim_O d(x,y)$.
Since $d(\mathfrak{X},\mathcal{Y}) \leq \delta$, there are
$x=x'_0,x'_1, \ldots , x'_{n-1},x'_n=y$ in $\mathcal{Y}$ such that
$d(x_i',x_i) \leq 3 / 2 \delta$ for $i=1, \ldots n-1$. Therefore, for
each $i$,
$d(x'_i,x'_{i+1}) \leq d(x'_i, x_i) + d(x_i, x_{i+1}) + d(x_{i+1},
x'_{i+1}) \leq 4 \delta$.~$\blacksquare$

The $\delta$-graph associated to a metric space $(\mathfrak{X},d)$ is
defined as follows; The vertices of the graph consists of the points
of $\mathfrak{X}$ and edges $x \longleftrightarrow y$ iff $d(x,y) \leq
\delta$.  Do not confuse the $\delta$-graph distance $\opr{dist}$ with
$\delta \times \opr{dist}$!

\begin{prop} \label{graph-metric-equiv-prop}
Suppose $(\mathfrak{X},d)$ is $\delta$-connected and $\opr{dist}$ is the
$\delta$-graph distance on $\mathfrak{X}$. Then there is a $c \gg 0$ such that
for all $x, y \in \mathfrak{X}$ with $\delta \leq d(x,y)$,
$$c \, \delta \, \opr{dist}(x,y) \leq d(x,y) \leq \delta \,
\opr{dist}(x,y).$$
\end{prop}
{\sc Proof.}  If $\opr{dist}(x,y) = n$, then by the definition of
graph distance, there is a sequence $x=x_0, \ldots , x_n=y$ such that
$d(x_i, x_{i+1}) \leq \delta$. Thus by the triangle inequality,
$d(x,y) \leq n \delta = \delta \, \opr{dist}(x,y)$.  In the other
direction, suppose $\delta \leq d(x,y)$.  $\opr{dist}(x,y)$ is the
smallest hyperinteger $n$ such that there is a sequence $x=x_0, x_1,
\ldots , x_n=y$ in $\mathfrak{X}$ with $d(x_i, x_{i+1}) \leq \delta$
for all $0 \leq i < n-1$.  Since $(\mathfrak{X},d)$ is a
$\delta$-space, there is a sequence with these properties and the
additional property $0 \ll d(x,y)/(n \, \delta)$.  In particular,
$c_{x,y} = d(x,y)/(\opr{dist}(x,y) \, \delta) \gg 0$.  Take $c = \inf
\{c_{x,y} : \delta \leq d(x,y)\} \gg 0$ which is the infimum over an
internal condition. Thus $d(x,y) \geq c \, \delta \,
\opr{dist}(x,y)$.~$\blacksquare$

\begin{prop} \label{graph-equivalent-to-original}
Let $(\mathfrak{X},d)$ be internal and $\delta$-connected. Then there
is a $\delta$-separated $\mathcal{Y} \subseteq \mathfrak{X}$ with
$d(\mathfrak{X},\mathcal{Y}) \leq \delta$ satisfying the following
property: If $\opr{dist}$ is the $\delta$-graph metric on
$\mathcal{Y}$, then $\delta \opr{dist}$ is $S$-Lipschitz equivalent to
$d$ on $\mathcal{Y}$.
\end{prop}
{\sc Proof.} Let $\mathcal{Y} \subseteq \mathfrak{X}$ be maximal
$\delta$-separated. In particular, $d(\mathfrak{X}, \mathcal{Y}) \leq
\delta$.  By Proposition~\ref{delta-subsets-prop}, $\mathcal{Y}$ is $4
\delta$-connected. Now if $x,y \in \mathcal{Y}$ with $x \neq y$, then
$d(x,y) \geq \delta$ so by Proposition~\ref{graph-metric-equiv-prop},
$$ 0 \ll \frac{d(x,y)} {4 \delta \opr{dist}(x,y)} \sim_O \frac{d(x,y)}
{\delta \opr{dist}(x,y)} \ll \infty.$$
Thus $\delta \, \opr{dist}$ on $\mathcal{Y}$ is $S$-Lipschitz
equivalent to the metric $d$.~$\blacksquare$

Note that nothing said so far implies $\mathcal{Y}$ can be taken to be
hyperfinite.

\begin{prop} \label{existence-of-approximating-graph}
Suppose $\delta \cong 0$ and $(\mathfrak{X},d)$ is an internal
$\delta$-connected space such that all closed balls
$\opr{\overline{B}}(x,\delta)$ are $\transfer$-precompact.  Then there
is a graph $\mathcal{Y} \subseteq \mathfrak{X}$ with
$d(\mathfrak{X},\mathcal{Y}) \leq \delta$ such that $\delta
\opr{dist}$ is $S$-Lipschitz equivalent to $d$ and $\mathcal{Y}$ is
hyperfinitely branching.  In particular $\opr{B}(a,r) \cap
\mathcal{Y}$ is hyperfinite for any $r \in \transfer \mathbb{R}$ and
$a \in \mathcal{Y}$.
\end{prop}
{\sc Proof.}  Let $\mathcal{Y} \subseteq \mathfrak{X}$ be as in
Proposition~\ref{graph-equivalent-to-original} with the $\delta$-graph
structure.  All nodes $a \in \mathcal{Y}$ have hyperfinitely many
adjacent nodes.  In fact let $A \subseteq
\opr{\overline{B}}(a,\delta)$ be hyperfinite and such that every $x
\in \opr{\overline{B}}(a,\delta)$ is distance $< \delta/2$ from some
$x_A \in A$. Such a set $A$ exists by $\transfer$-precompactness of
$\opr{\overline{B}}(a,\delta)$.  By definition of the $\delta$-graph,
all nodes adjacent to $a$ in the $\delta$-graph of $\mathcal{Y}$ are
members of $\opr{\overline{B}}(a,\delta)$; since $\mathcal{Y}$ is
$\delta$-separated, all nodes in $\mathcal{Y}$ are at distance $\geq
\delta$ from each other. Thus $x \mapsto x_A$ is injective on the set
of nodes of $\mathcal{Y}$ adjacent to $a$.  In particular the
cardinality of this set is $\leq \opr{card} A$ which is hyperfinite.
Therefore $\mathcal{Y}$ is hyperfinitely branching.  It follows that
the set of points of $\mathcal{Y}$ at graph distance $< \infty$ (but
possibly unlimited) from any $a \in \mathcal{Y}$ is hyperfinite.  In
particular, the set of points at $d$ distance $< \infty$ from any $a$
is hyperfinite.~$\blacksquare$

\paragraph{Remarks.} The graph $\mathcal{Y}$ may have unlimited
branching at each node. If $\mathfrak{X}$ has the property that
$\opr{\overline{B}}(x,\delta)$ is covered by a limited number of balls
$\opr{B}(x,\delta/2)$, then we may assume the graph $\mathcal{Y}$ has
limited branching.  For the set $A$ in the above proof can be assumed
limited.

\begin{thm}\label{main-characterization-result}
Suppose $(X,d)$ is an Ahlfors regular length space of dimension
$\lambda$ up to $R \in ]0, \infty]$.  Then there is a hyperfinite
graph space $(\mathfrak{X}, d_{\mathfrak{X}})$, a Lipschitz equivalent
metric $d'$ on $X$ and a hyperreal $R_{+}$ such that the following
hold:
\begin{enumerate}
\item $(X,d')$ is a bounded component of $\lhull{(\mathfrak{X},
d_{\mathfrak{X}})}$.
\item For some uniform hypermeasure $\nu$ on $\mathfrak{X}$,
$\opr{P}_\mathfrak{X}(\nu) |_X = \mathcal{H}^\lambda$ (Hausdorff
measure with respect to the original metric $d$, though see the
remarks below).
\item $(\mathfrak{X}, d_{\mathfrak{X}})$ is of polynomial growth
$\lambda$ on $[R_{-}, R_{+}]$ where $R_{-} \cong 0$ and $0 \ll R_{+} \ll
\infty$ if $R$ is finite or $R_{+}$ is an unlimited hyperreal in case $R
= \infty$.
\end{enumerate}
\end{thm}
{\sc Proof.}  Since $(X,d)$ is a complete Ahlfors regular space, it is
locally compact (See Proposition~\ref{ahlfors-condition}).  By basic
results on length spaces, all closed balls in $(X,d)$ are compact.  It
follows $\transfer (X,d)$ is an internal length space such that all
closed balls are $\transfer$-compact.  Moreover, by
Proposition~\ref{characterization-of-self-bounded-component}, the
bounded component of $X$ in $\phull{\transfer X}$ is $X$ itself.

Let $\delta$ be an arbitrary positive infinitesimal.  Consider the
$\delta$-graph structure on $\transfer (X,d)$.  Apply
Proposition~\ref{existence-of-approximating-graph} to $\mathfrak{X} =
\transfer X$ to conclude the existence of a subgraph $\mathcal{Y}$
with the properties stated there.  In particular, the metrics $\delta
\, \opr{dist}$ and $\transfer d$ on $\mathcal{Y}$ are $\ESS$-Lipschitz
equivalent. Since $d_{\transfer X}(\transfer X, \mathcal{Y}) \leq
\delta$, if follows $\lhull{\mathcal{Y}} = \phull{\transfer X}$.
Moreover, the metric $\lhull{(\delta \opr{dist})}$ is Lipschitz
equivalent to $\phull{\transfer d}$.  Now let $a \in \mathcal{Y}$ be
such that $\hat{a} \in X$; then for $r$ unlimited (but $r \in
\transfer \mathbb{R}$), $\mathcal{Z} = \opr{B}(a, r) \cap \mathcal{Y}$
is hyperfinite, and such that $\lhull{\mathcal{Z}} \supseteq X$ and
$X$ is a bounded component of $\lhull{\mathcal{Z}}$.

Let $d'$ be the metric $\lhull{(\delta \opr{dist})}$ restricted to
$X$. By the previous paragraph, $d'$ is Lipschitz equivalent to
$\phull{\transfer d} |_X = d$.  Now $(X,d)$ is separable since it is a
complete length space and in particular, $X$ is $\sigma$-compact with
respect to $d'$. Let $\mu$ be $\lambda$-dimensional Hausdorff measure
relative to the metric $d$.  $\mu$ is finite on compact sets and is
boundedly equivalent to Hausdorff measure relative to the metric $d'$.
It follows by Proposition~\ref{measlift-unbounded}
that there is a hypermeasure $\nu$ on $\mathcal{Z}$ such that
$\opr{P}_\mathcal{Z}(\nu)|_X = \mu$.  Now $\opr{supp}(\nu) = \{z \in
\mathcal{Z}: \nu_z \neq 0 \}$ is an internal set.  In addition, by the
Ahlfors regularity property of Hausdorff measure,
$\phull{\opr{supp}(\nu)} \supseteq X$.  Thus we can assume without
loss of generality that $\nu_z \neq 0$ for all $z \in \mathcal{Z}$.
By Proposition~\ref{even-covering-metric-case}, we can also assume
$\nu$ is a uniform hypermeasure on $\mathcal{Z}$ with density $1/M$.
%
%
Applying Proposition~\ref{ahlfors-regular-characterization} and the
remarks after the proof of
Proposition~\ref{ahlfors-regular-characterization} to $(\mathcal{Z},
d_{\mathcal{Z}})$, we can immediately read off the main conclusion of
the theorem.  Thus there is a hyperfinite $\mathcal{Z}' \subseteq
\mathcal{Z}$ such that $\lhull{\mathcal{Z}'} \supseteq X$ and
$(\mathcal{Z}',d_{\mathcal{Z}})$ is of polynomial growth on some
interval $[R_{-}, R]$ where $R_{-} \cong 0$.~$\blacksquare$

\paragraph{Remarks.}  Note that in the previous result,
$\mathcal{H}^\lambda$ may be replaced by any measure which is
boundedly equivalent to $\mathcal{H}^\lambda$.  See the remark after
Proposition~\ref{ahlfors-regular-characterization}.

\bibliography{../../../nsa.dir/refs}
\bibliographystyle{plain}

\noindent 9112 Decatur Ave S. \\ Bloomington, MN 55438  \\
{\em E-mail address:}  \tt{jt@mitre.org}

\end{document}